\documentclass[12pt,oneside]{article}

\usepackage{cmap}

\usepackage[cp1251]{inputenc}
\usepackage{amsthm,amsfonts,amsbsy, amssymb,amsmath,graphicx,euscript}
\usepackage{makeidx}
\usepackage{multicol}

\usepackage{mathrsfs}

\usepackage[unicode]{hyperref}

\usepackage{amsfonts}
\usepackage{cmap}
\usepackage{amssymb}
\usepackage{amsmath}
\usepackage{epsfig}
\usepackage{psfrag}
\usepackage{floatflt}
\usepackage{mathrsfs}
\usepackage{indentfirst}
\usepackage{cmap}
\author{И. Д. Рeмизoв}

\newcommand{\tr}{\mathrm{tr}}

\newcommand{\R}{\mathbb{R}}
\newcommand{\N}{\mathbb{N}}

\newcommand{\be}{\begin{equation}}
\newcommand{\ee}{\end{equation}}

\newcommand{\QED}{ \begin{flushright}$\Box$ \end{flushright}}

\theoremstyle{definition}
\newtheorem{df}{Definition}[section]
\newtheorem{pr}{Proposition}[section]

\newtheorem{rem}{Remark}[section]

\newtheorem{lm}{Lemma}[section]
\newtheorem{thm}{Theorem}[section]

\usepackage{ulem}
\usepackage{soul}

\usepackage{geometry}
\begin{document}


\begin{center}
\textbf{Solution to a parabolic differential equation in Hilbert space\\ via Feynman formula  - parts I and II}
\end{center}

\begin{center}
I.\,D.~Remizov\footnote{Email: ivremizov@yandex.ru}
\end{center}

\footnotesize
A parabolic partial differential equation $u'_t(t,x)=Lu(t,x)$ is considered, where $L$ is a linear second-order differential operator with time-independent coefficients, which may depend on $x$. We assume that the spatial coordinate $x$ belongs to a finite- or infinite-dimensional real separable Hilbert space $H$. 

Assuming the existence of a strongly continuous resolving semigroup for this equation, we construct a representation of this semigroup by a Feynman formula, i.e. we write it in the form of the limit of a multiple integral over $H$ as the multiplicity of the integral tends to infinity. This representation gives a unique solution to the Cauchy problem in the uniform closure of the set of smooth cylindrical functions on $H$. Moreover, this solution depends continuously on the initial condition. In the case where the coefficient of the first-derivative term in $L$ vanishes we prove that the strongly continuous resolving semigroup exists (this implies the existence of the unique solution to the Cauchy problem in the class mentioned above) and that the solution to the Cauchy problem depends continuously on the coefficients of the equation.\\ 

MSC2010 codes: 35C15, 47D06,	28C20. UDC code: 517.987.4.
\normalsize

\tableofcontents

\section{Introduction}

Representation of a function by the limit of a multiple integral as multiplicity tends to infinity is called a Feynman formula,
after the inventor of the equations of such type, R.P. Feynman, who was the first to use them (on the physical level of rigor)
for the solution of the Cauchy problem for PDEs \cite{F1, F2}. The term "Feynman formula" in this sense was introduced in 2002
by O.G. Smolyanov \cite{Sm2}. One can find out more about the Feynman formulas' research up to 2009 in \cite{Sm3}. It is important
to note that Feynman formulas are closely related to Feynman-Kac formulas \cite{Simon}, however the latter will not be studied in the present article. Usage of Feynman and Feynman-Kac formulas includes exact or numerical  evaluation of integrals over Gaussian measures on spaces of high or infinite dimension; some useful approaches to this topic are developed in \cite{EgZhLob, Lobanov}.

Differential equations for functions of an infinite-dimensional argument arise in (quantum) field theory and string theory,
theory of stochastic processes and financial mathematics.
Evolutionary equations (i.e. PDEs in the form $u'_t(t,x)=\dots$) in infinite-dimensional spaces have been studied since 1960s
by O.G. Smolyanov, E.T. Shavgulidze, E. Nelson, A.Yu. Khrennikov, S. Albeverio, L.C.L. Botelho and others. We will mention just some of the publications, which
are most recent and relevant for our study. 

In \cite{Butko1} the Schr\"{o}dinger equation in Hilbert space is studied. The equation includes the terms of second, first and zero order,
the coefficient of the second order term is constant. The solution to the Cauchy problem is given by a Feynman-Kac-Ito formula.

In \cite{DaPrato2} a solution to a heat equation in Hilbert space without the terms of the first and zero order is discussed,
the coefficient of the second-derivative term is constant. The solution is given in the form of a convolution with
the Gaussian measure (analogous to the finite dimensional equation with constant coefficients),
the existence of the resolving semigroup is proved. In \cite{Botelho} the solution to the same equation is given by a Feynman-Kac formula.

In \cite{BGS2010} the parabolic equation in finite-dimensional space is studied for the case of variable coefficients.
Under the assumption that a strongly continuous resolving semigroup exists for the Cauchy problem, Feynman and Feynman-Kac formulas were proven
in \cite{BGS2010} for the solution.

In \cite{Remizov}, for a class of equations in an infinite-dimensional space, with a variable coefficient at the highest derivative (but without first- and zero-order derivatives' terms),
a Feynman formula was obtained and the existence of resolving semigroup was proven.

In spaces over the field of p-adic numbers, Feynman and Feynman-Kac formulas for the solutions of the Cauchy problem
for evolutionary equations were given in \cite{SShK, SSham}.

In \cite{Plya1, Plya2}, Schr\"{o}dinger and heat equations in $\R^n$ were studied
in the case of time-dependent coefficients, and a Chernoff-type theorem was proven for this case.

In \cite{Butko-mult, BShS} Feynman formulas for perturbed semigroups are obtained.

In \cite{Bot1, Bot2} non-linear and semilinear heat-type equations are studied with the path integral approach, and in \cite{Bot3,Bot4} this approach is applied to the wave equation.\\

The present article extends my first results in this area \cite{Remizov} to the case of non-zero coefficients at the first- and zero-order derivatives.

\section{Notation and definitions}\label{denotions}

The symbol $H$ stands for the real separable Hilbert space with the scalar product $\left<\cdot,\cdot\right>$.

The self-adjoint, positive, non-degenerate (hence injective), linear operator $A\colon H\to H$ is assumed to be defined everywhere on $H$.
The operator $A$ is assumed to be of trace class, which means that for every orthonormal basis  $(e_k)$ in $H$ the sum
$\sum_{k=1}^\infty\left<Ae_k,e_k\right>=\tr A$ is finite; this sum is called the trace of $A$ (it is independent of
the choice of the basis $(e_k)$).

The symbol $\mathcal{X}$ below stands for any complex Banach space.
The symbol $L_b(\mathcal{X},\mathcal{X})$ stands for space of all linear bounded operators in $\mathcal{X}$, endowed with the
classical operator norm.

Symbol $C(M,N)$ will mean the set of all continuous functions from $M$ to $N$, where $M$ and $N$ are topological spaces.

A function $f\colon H\to \R$ is called cylindrical \cite{DF, SSh}, if there exist vectors $e_1,\dots,e_n$ from $H$ and function $f^n\colon \R^n\to\R$ such
that for every $x\in H$ the equality
$f(x)=f^n(\left<x,e_1\right>,\dots, \left<x,e_n\right>)$ holds. In other words,
the function $f\colon H\to  \R$ is cylindrical if there exists an $n-$dimensional subspace
$H_n\subset H$ and orthogonal projector $P\colon H\to H_n$ such that
$f(x)=f(Px)$ for every $x\in H$. The cylindrical function $f$
can be imagined as a function, which is first defined on $H_n$ and then continued to the entire space $H$ in such a way that $f(x)=f(x_0)$ if
$x_0\in H_n$ and $x\in (x_0+\mathrm{ker}P)$.

Symbol $D=C^\infty_{b,c}(H,\R)$ stands for the space of all continuous bounded cylindrical functions $H\to \R$ such that they
have Fr\'{e}chet derivatives \cite{Cartan} of all positive integer orders at every point of  $H$, and their Fr\'{e}chet derivatives of any positive integer order are bounded and continuous.

If $f\colon H\to\R$ is twice Fr\'{e}chet differentiable, then $f'(x)$ will stand for the first
Fr\'{e}chet derivative of $f$ at the point $x$, and $f''(x)$ will denote the second derivative.
Riesz-Fr\'{e}chet representation theorem allows us to assume $f'(x)\in H$ and $f''(x)\in L_b(H,H)$ for every $x\in H$.

Symbol $C_b(H,\R)$ stands for the Banach space of all bounded continuous functions $H\to \R$,
endowed with a uniform norm $\|f\|=\sup_{x\in H}|f(x)|$. It is regarded as a closed subspace of a complex Banach
space $C_b(H,\mathbb{C})$.

Let $X=\overline{C^\infty_{b,c}(H,\R)}$ be the closure of the space $D$ in $C_b(H,\R)$. It is clear, that $X$
with the norm $\|f\|=\sup_{x\in H}|f(x)|$ is a Banach space, as it is a closed linear subspace of the Banach space $C_b(H,\R)$.
Function $f$ belongs to $X$ if and only if there is a sequence of functions $(f_j)\subset D$ such that $\lim_{j\to\infty}f_j=f$, i.e.
$\lim_{j\to\infty}\sup_{x\in H}|f(x)-f_j(x)|=0$.

Symbol $C_b(H,H)$ stands for a Banach space of all bounded continuous functions $B\colon H\to H$,
endowed with the uniform norm $\|B\|=\sup_{x\in H}\|B(x)\|$.

Denote $D_H=\{B\colon H\to H \big |  \exists N\in\N , b_k\in H, B_k\in D: B(x)=B_1(x)b_1+\dots+B_N(x)b_N\}$.

Let $X_H$ be the closure of $D_H$ in $C_b(H,H)$.

If $x\in H$, and $R\colon H\to H$ is linear, trace class, positive, non-degenerate operator, then symbol $\mu_R^x$ stands for
the Gaussian probabilistic measure \cite{Bog2, DF, Go} on $H$ with expectation $x$ and correlation operator $R$, i.e. the unique sigma-additive
measure on Borel sigma-algebra in $H$ such that the
equality $\int_H e^{i\left<z,y\right>}\mu^x_R(dy)=\exp\left(i\left<z,x\right>-\frac{1}{2}\left<Rz,z\right>\right)$ holds for every $z\in H$. To make it shorter,
we will write $\mu_R$ instead of $\mu_R^0$.

If $B\colon H\to H$ is a vector field, and $g\colon H\to\R$ and $C\colon H\to\R$ are real-valued functions, then symbol $L$ defines a differential
operator on the space of functions $\varphi\colon H\to \R$

$$ (L\varphi)(x):=g(x)\tr A\varphi''(x)+\langle \varphi'(x),AB(x)\rangle + C(x)\varphi(x),\quad x\in H.$$

The pair $(\mathcal{L},M)$ defines a linear operator $\mathcal{L}$ with the domain $M$. It will be shown in theorem \ref{gDeltaA_thm} that
$L(D)\subset X$ when $A$, $B$, $g$ and $C$ have certain properties. So $(L, D)$ is a densely
defined (on $D$) operator $L\colon X\supset D \to X$. Here the earlier defined spaces $D$ and $X$ are endowed with the uniform norm,
induced from $C_b(H,\R)$. Let $(\overline{L}, D_1)$ be the closure of $(L, D)$ in $X$. This means that
$$D_1=\{f\in X\big|\exists (f_j)\subset D: \lim_{j\to\infty}f_j=f, \exists\lim_{j\to\infty}Lf_j\},$$
and, if $f\in D_1$, then, by definition, $\overline{L}f=\lim_{j\to\infty}Lf_j.$

If for every fixed first argument $t>0$ of the function $u\colon [0,+\infty)\times H\to\R$ we have
$[x\longmapsto u(t,x)]\in D_1$, then the expression $\overline{L}u(t,x)$ means the result of applying
the operator $\overline{L}$ to the function $x\longmapsto u(t,x)$ with the fixed $t> 0$.

Expression $(S_t)_{t\geq 0}$ defines the one-parameter family of linear operators in the space of functions $\varphi\colon H\to \R$
$$ (S_t\varphi)(x):= e^{tC(x)}\int_H\varphi(x+y)e^{\left<\frac{1}{g(x)}B(x),y\right>}\mu_{2tg(x)A}(dy) \textrm{ for }t>0, \textrm{ and }S_0\varphi:=\varphi.$$

\rem Further, in theorem
\ref{preChern2}, we will prove that for every $t\geq 0$ and for $A$, $B$, $g$ and $C$
having certain properties the following holds i) $S_t(X)\subset X$,
ii) operator $S_t$ is bounded, and iii) $\frac{d}{dt}S_t\varphi\big|_{t=0}=L\varphi$ for all $\varphi\in D$.
This will allow us to use the Chernoff approximation (theorems \ref{FormulaChernova}, \ref{FormulaChernova2})
and prove the main result of the present article, theorem \ref{main_th_art}.

\section{Helpful facts and techniques}

\subsection{Integration in Hilbert space}

\lm\label{integrir_cil_f} (\cite{DF}, Chapter II, \S 2, $3^\circ$) If a function $\varphi\colon H\to\mathbb{R}$ is cylindrical and measurable,
i.e. $\varphi(x)=\varphi^n(\langle x,e_1\rangle,\dots,\langle x,e_n\rangle)$ for some $n\in\mathbb{N}$,
some measurable function $\varphi^n\colon\mathbb{R}^n\to\mathbb{R}$,
and some finite orthonormal family of vectors $e_1,\dots,e_n$ from space $H$, then
\be\label{formula_odin} \int_H \varphi(y)\mu_A(dy)=\left(\frac{1}{\sqrt{2\pi}}\right)^n\frac{1}{\sqrt{\det M_Q}}\int_{\mathbb{R}^n}
\varphi^n(z) \exp\left(-\frac{1}{2}\left<M_Q^{-1}z,z\right>_{\mathbb{R}^n}\right)dz, \ee
where $H_n=\mathrm{span}( e_1,\dots, e_n)$,  and
$P\colon H\ni h\longmapsto \langle h,e_1\rangle e_1+\dots+\langle h,e_n\rangle e_n\in H_n,$
$Q=PA$, $Q\colon H_n\to H_n$, and $M_Q$ is the matrix of the operator $Q$ in basis
$e_1,\dots,e_n$ of the space $H_n$.
If
$e_1,\dots,e_n$
is a full set of eigenvectors of the operator $Q$, and $q_1,\dots,q_n$ is the corresponding set of eigenvalues, then
\be\label{intcyl} \int_H \varphi(y)\mu_A(dy)=\left(\frac{1}{\sqrt{2\pi}}\right)^n\frac{1}{\sqrt{\prod_{i=1}^n q_i}}
\int_{\mathbb{R}^n} \varphi^n(z_1,\dots, z_n) \exp\left(-\sum_{i=1}^n \frac{z_i^2}{2q_i} \right)dz_1\dots dz_n. \ee

\lm (\textit{Explicit form of some integrals over Gaussian measure})

Let $H$ be a real separable Hilbert space of finite or infinite dimension,
$\widetilde{A}\colon H\to H$ be a linear, trace class, symmetric, positive, non-degenerate operator,
$\mu_{\widetilde{A}}$ be the centered Gaussian measure on $H$ with the correlation operator $\widetilde{A}$, and $G\colon H\to H$
be a bounded linear operator. Let $w$ and $z$ be non-zero vectors from $H$.

Then the following equalities hold:

\be\label{int_kv_form1} \int_H\langle Gy,y\rangle \mu_{\widetilde{A}}(dy)=\tr(\widetilde{A}G), \ee

\be\label{int_eksp} \int_H e^{\langle z,y\rangle} \mu_{\widetilde{A}}(dy)=e^{\frac{1}{2}\langle \widetilde{A}z,z\rangle}, \ee

\be\label{int_eksp1}\int_H \langle w, y\rangle e^{\langle z,y\rangle} \mu_{\widetilde{A}}(dy)= \langle \widetilde{A}w,z\rangle e^{\frac{1}{2}\langle \widetilde{A}z,z\rangle}, \ee

\be\label{int_eksp2}\int_H \langle Gy, y\rangle e^{\langle z,y\rangle} \mu_{\widetilde{A}}(dy)= (\tr \widetilde{A}G + \langle G\widetilde{A}z,\widetilde{A}z\rangle)
 e^{\frac{1}{2}\langle \widetilde{A}z,z\rangle}. \ee

\textbf{Proof.} Formulas (\ref{int_kv_form1}) and (\ref{int_eksp}) can be found in \cite{DF}, chapter II, \S 2, $1^\circ$.
Formula (\ref{int_eksp1}) can be derived from the fact that the function under the integral is cylindrical, so
lemma \ref{integrir_cil_f} can be employed. For a proof of (\ref{int_eksp2}), one can make the change of variable in the integral, $h=y-Aw$, then
(\cite{DF}, chapter II, \S 4, $2^\circ$, theorem 4.2) we have $\mu_{\widetilde{A}}(dy)=e^{-\frac{1}{2}\left<\widetilde{A}w,w\right>-\left<h,w\right>}\mu_{\widetilde{A}}(dh)$,
and the integral reduces to (\ref{int_kv_form1}).

\lm\label{linzam} (\textit{On a linear change of variable in the integral over Gaussian measure})
Let $H$ be a real separable Hilbert space.
Suppose a linear operator $A\colon H\to H$ is positive, non-degenerate, trace class, and self-adjoint.
We will identify with the symbol $\mu_A$ the centered Gaussian measure on $H$ with the correlational operator $A$.
Let $t>0$; the symbol $tA$ denotes operator, that takes $x\in H$ to $tAx\in H$.
Let $f\colon H\to\mathbb{R}$ be a continuous integrable function.

Then \be\label{lz}\int_Hf(x)\mu_{tA}(dx)=\int_Hf(\sqrt{t}x)\mu_A(dx).\ee

\textbf{Proof} uses the uniqueness of the Gaussian measure with a given Fourier transform, and the standard theorem of changing variable in the Lebesgue integral.

\lm\label{polyexpint}(\textit{On integrability of a polynomial multiplied by an exponent})
Let $H, A, \mu_A$ be as above, $P\colon \R\to\R$ be a polynomial, and $\beta\in\R$.

Then function $H\ni x\longmapsto P(\|x\|)e^{\beta\|x\|}\in\R$ is integrable over $\mu_A$.

\textbf{Proof} is easy to construct by relying on Fernique's theorem \cite{Fern}, which (applied to this case)
says that there exists such $\alpha>0$ that $\int_H e^{\alpha\|y\|^2}\mu_A(dy)<+\infty$.

\subsection{Differentiation in Hilbert space}

\pr\label{pr_diff} Let $f$ be a cylindrical real-valued function on $H$, i.e. there is a number $n\in \N$ and
a function $f^n\colon \R^n\to\R$ such that for every $x\in H$
the equality $f(x)=f^n(\langle x,e_1\rangle,\dots,\langle x,e_n\rangle)$ holds.
A set of vectors $e_1,\dots, e_n$ can be considered orthonormal without loss of generality.
Lets complete this set to an orthonormal basis $(e_k)_{k\in \N}$ in $H$.

Then:

1. Function $f$ is differentiable in the direction $h$ if and only if the function $f^n$ is differentiable in the direction
$(\langle h,e_1\rangle,\dots,\langle h,e_n\rangle)\in\R^n$, and
$$ f'(x)h=\left\langle h, \Big(\partial_1f^n(\langle x,e_1\rangle,\dots,\langle x,e_n\rangle),\ \dots,\ \partial_nf^n(\langle x,e_1\rangle,\dots,\langle x,e_n\rangle),\ 0,\ 0,\ 0,\dots\Big)  \right\rangle,$$
where the symbol $\partial_j f^n$ defines the partial derivative with respect to the $j$-th argument of the function
$f^n$, and $(\alpha_1,\dots,\alpha_n,0,0,0,\dots)=\alpha_1e_1+\dots+\alpha_ne_n$.
If the function $f$ has a Fr\'{e}chet derivative at the point $x$, then $f'(x)$
is a vector whose first $n$ coordinates yield the gradient of the function $f^n$,
and the other coordinates are zero:
\be\label{proizv_napr}f'(x)=
\Big(\partial_1f^n(\langle x,e_1\rangle,\dots,\langle x,e_n\rangle),\ \dots,\ \partial_nf^n(\langle x,e_1\rangle,\dots,\langle x,e_n\rangle),\ 0,\ 0,\ 0,\dots\Big). \ee

2. Function $f$ has a Fr\'{e}chet derivative in $H$ if and only if the function $f^n$
has a Fr\'{e}chet derivative in $\R^n$.

3. Let $A\colon H\to H$ be a trace-class operator (i.e. let $\tr A<\infty$). Then
\be\label{traf}\tr Af''(x)=\sum_{s=1}^n\sum_{k=1}^n \langle Ae_s,e_k\rangle
\Big(\partial_k\partial_s f^n(\langle x,e_1\rangle,\dots,\langle x,e_n\rangle)\Big)=\tr \Big(A_n(f^n)''(\langle x,e_1\rangle,\dots,\langle x,e_n\rangle)\Big),\ee

where $A_n$ is the matrix of the operator $PA$ in the basis $e_1,\dots,e_n$, where $P$ is the projector to the linear span of the vectors $e_1,\dots,e_n$.

\textbf{Proof} is a straight-forward application of the derivative's definition.

\pr For $(n+1)$-times Fr\'{e}chet differentiable function $f\colon H\to \mathbb{R}$ there is \cite{Sm1} a Taylor decomposition

\be\label{Teyl} f(x+h)=f(x)+\frac{1}{1!}f'(x)h+\frac{1}{2!}f''(x)(h,h)+\dots+\frac{1}{n!}f^{(n)}(x)(h,\dots,h)+R_n(x,h),\ee

and  \be\label{oc_ost_chl}|R_n(x,h)|\leq\frac{\|h\|^{n+1}}{(n+1)!}\sup\limits_{z\in[x,x+h]}\left\|f^{(n+1)}(z)\right\|.\ee

\subsection{Differential operator on a finite-dimensional space}

\lm\label{lm_urchp} (\cite{Kryl}, theorems 4.3.1, 4.3.2. and Corollary 4.3.4) Suppose for every $i=1,\dots,n$ and $j=1,\dots,n$
functions $a^{ij}\colon\R^n\to \R$, $b^{i}\colon\R^n\to \R$, $c\colon\R^n\to \R$ from $C^\infty_b(\R^n,\R)$ are given, where
$C^\infty_b(\R^n,\R)$ is the class of all bounded real-valued functions on $\R^n$, which have bounded partial derivatives of all orders.
Suppose also that $c(x)\leq 0$ for all $x\in\R^n$.

For $u\in C^\infty_b(\R^n,\R)$ we define a differential operator $T$ by the formula

$$(Tu)(x)=\sum_{i=1}^n\sum_{j=1}^na^{ij}(x) \frac{\partial^2}{\partial x_i\partial x_j}u(x) + \sum_{i=1}^nb^i(x)\frac{\partial}{\partial x_i}u(x)+c(x)u(x).$$

Suppose that there exists a constant $\varkappa>0$ such that for every $\xi=(\xi_1,\dots,\xi_n)\in\R^n$ and all $x\in\R^n$
the ellipticity condition is fulfilled:  $\sum_{i=1}^n\sum_{j=1}^na^{ij}(x)\xi_i\xi_j\geq\varkappa\|\xi\|^2.$
Take an arbitrary constant $\lambda>0$ and function $f\in C^\infty_b(\R^n,\R)$.

Then:

1. There is a unique function $u\in C^\infty_b(\R^n,\R)$, which is a solution of the equation
\be\label{urchp} (Tu)(x) - \lambda u(x)=f(x).\ee

2. For every function $v\in C^\infty_b(\R^n,\R)$ the following estimate is true

\be\label{diss_oc} \sup_{x\in\R^n} \left|(Tv)(x) - \lambda v(x)\right|\geq \lambda\sup_{x\in\R^n}|v(x)|.\ee

Note that equation (\ref{urchp}) can have unbounded solutions; this does not contradict the lemma.

\subsection{Strongly continuous semigroups of operators and evolutionary equations}

Let $\mathcal{X}$ be a complex Banach space.

\df By a $C_0$-semigroup, or a strongly continuous one-parameter semigroup $(T_s)_{s\geq 0}$ of linear bounded operators in $\mathcal{X}$ we (following \cite{EN, Pazy}) mean the mapping
$$T\colon [0,+\infty)\to L_b(\mathcal{X},\mathcal{X})$$
of the non-negative  half-line into the space of all bounded linear operators on $\mathcal{X}$,
which satisfies the following conditions:

1. $\forall \varphi\in \mathcal{X}: T_0\varphi=\varphi.$

2. $\forall t\geq 0,\forall s\geq 0: T_{t+s}=T_t\circ T_s.$

3. $\forall \varphi\in \mathcal{X}$ function $s\longmapsto T_s\varphi$ is continuous as a mapping $[0,+\infty)\to \mathcal{X}.$

\df By the generator of a strongly continuous one-parameter semigroup $(T_s)_{s\geq 0}$ of linear bounded operators on $\mathcal{X}$ we mean
a linear operator $\overline{\mathcal{L}}\colon \mathcal{X}\supset Dom(\overline{\mathcal{L}})\to \mathcal{X}$
given by the formula

$$\overline{\mathcal{L}}\varphi = \lim_{s\to +0}\frac{T_s\varphi-\varphi}{s}$$

on its domain
$$Dom(\overline{\mathcal{L}})=\left\{\varphi\in \mathcal{X}: \exists \lim_{s\to +0}\frac{T_s\varphi-\varphi}{s}\right\},$$

where the limit is understood in the strong sense, i.e. it is defined in terms of the norm in the space $\mathcal{X}$.\\

The use of the symbol $\overline{\mathcal{L}}$ for the generator is related to the fact that the generator is always a closed operator:

\pr (theorem 1.4 in \cite{EN}, p. 51) The generator of a strongly continuous semigroup is a closed linear operator with a dense domain.
The generator defines its semigroup uniquely.

\pr\label{path_diff} (lemma 1.1 and definition 1.2. in \cite{EN}, p. 48-49) The set $Dom(\overline{\mathcal{L}})$
coincides with the set of those  $\varphi\in \mathcal{X}$, for which the mapping $s\longmapsto T_s\varphi$ is differentiable with respect to  $s$
at every point  $s\in[0,+\infty).$

\df\label{abszadC} 1. The problem of finding a function $U\colon [0,+\infty)\to \mathcal{X}$ such that

\be\label{ACP1}
\left\{ \begin{array}{ll}
 \frac{d}{d t}U(t)= \overline{\mathcal{L}}U(t); & t\geq 0,\\
 U(0)=U_0,\\
  \end{array} \right.
\ee

is called the abstract Cauchy problem, associated
with the closed linear operator $\overline{\mathcal{L}}\colon \mathcal{X}\supset Dom(\overline{\mathcal{L}})\to \mathcal{X}$ and a vector $U_0\in \mathcal{X}$.

2. A function $U\colon [0,+\infty)\to \mathcal{X}$ is called a classic solution to abstract Cauchy problem (\ref{ACP1})
if, for every $t\geq 0$, the function $U$ has a continuous derivative $U'\colon [0,+\infty)\to \mathcal{X}$,
$U(t)\in Dom(\overline{\mathcal{L}})$, and (\ref{ACP1}) holds.

3. A continuous function $U\colon [0,+\infty)\to \mathcal{X}$ is called a mild solution
to abstract Cauchy problem (\ref{ACP1}) if  for every $t\geq 0$ we have $\int_{0}^tU(s)ds\in Dom(\overline{\mathcal{L}})$
and $U(t)=\overline{\mathcal{L}}\int_0^tU(s)ds + U_0$.

\pr\label{ACPsol} (proposition 6.2 in \cite{EN}, p. 145) If the operator $(\overline{\mathcal{L}}, Dom(\overline{\mathcal{L}}))$
is a generator of a strongly continuous semigroup $(T_s)_{s\geq 0}$, then:

1. For every $U_0\in Dom(\overline{\mathcal{L}})$ there is a unique classic solution to abstract Cauchy problem (\ref{ACP1}), which is
given by the formula $U(t)=T(t)U_0$.

2. For every $U_0\in \mathcal{X}$ there is a unique mild solution to abstract Cauchy problem (\ref{ACP1}), which is
given by the formula $U(t)=T(t)U_0$.

\df Linear operator $\mathcal{L}\colon \mathcal{X}\supset Dom(\mathcal{L})\to \mathcal{X}$ in Banach space $\mathcal{X}$ is called
dissipative if for every  $\lambda>0$ and every $x\in Dom(\mathcal{L})$
the estimate $\|\mathcal{L}x-\lambda x\|\geq \lambda\|x\|$ holds.

\pr\label{pr_dissop} (\textit{On the closability of a densely defined dissipative operator})
(proposition 3.14 in \cite{EN}) A linear dissipative operator $\mathcal{L}:\mathcal{X}\supset Dom(\mathcal{L})\to \mathcal{X}$
in the Banach space $\mathcal{X}$ with the domain $Dom(\mathcal{L})$ dense in $\mathcal{X}$ is closable.
The closure $\overline{\mathcal{L}}: \mathcal{X}\supset Dom(\overline{\mathcal{L}})\to X$ is also a dissipative operator.\\

The main tool for the construction of Feynman formulas for the solutions of the Cauchy problem is Chernoff's theorem. For convenience we decompose its conditions into several blocks and give them separate names, as follows.

\thm \label{FormulaChernova}(\textsc{P.\,R.~Chernoff, 1968}; see \cite{Chernoff} and theorem 10.7.21 in \cite{BS}) Let $\mathcal{X}$ be Banach space,
and $L_b(\mathcal{X}, \mathcal{X})$ be the space of all linear bounded operators in $\mathcal{X}$ endowed with the operator norm. Let $\overline{\mathcal{L}}\colon \mathcal{X}\supset Dom(\overline{\mathcal{L}})\to \mathcal{X}$ be a linear operator. 

\textbf{Suppose} there is a function $F$ such that:

(E). There exists a strongly continuous semigroup $(e^{t\overline{\mathcal{L}}})_{t\geq 0},$ and its generator is $(\overline{\mathcal{L}},Dom(\overline{\mathcal{L}}))$.

(CT1). $F$ is defined on $[0,+\infty)$, takes values in $L_b(\mathcal{X},\mathcal{X})$ and $t\longmapsto F(t)f$ is continuous for every vector $f\in\mathcal{X}$. 

(CT2). $F(0)=I$. 

(CT3). There exists a dense subspace $\mathcal{D}\subset \mathcal{X}$ such that for every $f\in \mathcal{D}$ there exists
a limit $F'(0)f=\lim_{t\to 0}(F(t)f-f)/t=\mathcal{L}f$. 

(CT4). The operator  $(\mathcal{L},\mathcal{D})$ has a closure $(\overline{\mathcal{L}},Dom(\overline{\mathcal{L}})).$ 

(N). There exists $\omega\in\mathbb{R}$ such that $\|F(t)\|\leq e^{\omega t}$ for all $t\geq 0$.

\textbf{Then} for every $f\in \mathcal{X}$ we have $(F(t/n))^nf\to e^{t\overline{\mathcal{L}}}f$ as $n\to \infty$, and the limit
is uniform with respect to $t$ from every segment $[0,t_0]$ for every fixed $t_0>0$.

\df\label{Cheq}In the present article two mappings $F_1$ and $F_2$ are called Chernoff-equivalent if there exists a $C_0$-semigroup $(e^{t\overline{\mathcal{L}}})_{t\geq 0}$ such that $(F_1(t/n))^nf\to e^{t\overline{\mathcal{L}}}f$, $(F_2(t/n))^nf\to e^{t\overline{\mathcal{L}}}f$ for every $f\in \mathcal{X}$ as $n\to \infty$, and the limit is uniform with respect to $t$ from every segment $[0,t_0]$ for every fixed $t_0>0$.

\rem There are several slightly different definitions of the Chernoff equivalence, see e.g. \cite{SWW, OSS, Butko4}. We will just use this one not going into details. The only thing we need from this definition is that if $F$ satisfies all the conditions of Chernoff's theorem, then by Chernoff's theorem the mapping $F$ is Chernoff-equivalent to the mapping $F_1(t)=e^{t\overline{\mathcal{L}}}$, i.e. the limit of $(F(t/n))^n$ as $n$ tends to infinity yields the $C_0$-semigroup $(e^{t\overline{\mathcal{L}}})_{t\geq 0}$.

\df\label{weakCher}Let us follow \cite{Rmz1} and call a mapping $F$ Chernoff-tangent to the operator $\mathcal{L}$ if it satisfies the conditions (CT1)-(CT4) of Chernoff's theorem.

\rem With these definitions the Chernoff-equivalence of $F$ to $(e^{t\overline{\mathcal{L}}})_{t\geq 0}$ follows from: existence (E) of the $C_0$-semigroup + Chernoff-tangency (CT) + growth of the norm bound (N).

\thm\label{FormulaChernova2}(Chernoff-type theorem, \cite{EN}, corollary 5.3 from theorem 5.2)
Let $\mathcal{X}$ be a Banach space, and $L_b(\mathcal{X}, \mathcal{X})$ be the space of all linear bounded operators on $\mathcal{X}$
endowed with the operator norm. Suppose there is a function
$$ V\colon [0,+\infty)\to L_b(\mathcal{X},\mathcal{X}),$$
meeting the condition $V_0=I$, where $I$ is the identity operator.
Suppose there are numbers $M\geq 1$ and $\omega\in\R$ such that $\|(V_t)^k\|\leq Me^{k\omega t}$ for every $t\geq 0$ and every $k\in \N$. Suppose the limit
$$\lim_{t\downarrow 0}\frac{V_t\varphi - \varphi}{t}=:\mathcal{L}\varphi$$
exists for every $\varphi\in \mathcal{D}\subset \mathcal{X},$ where $\mathcal{D}$ is a dense subspace of $\mathcal{X}$. Suppose there is 
a number $\lambda_0>\omega$ such that $(\lambda_0I-\mathcal{L})(\mathcal{D})$ is a dense subspace of $\mathcal{X}$.

Then the closure $\overline{\mathcal{L}}$ of the operator $\mathcal{L}$ is a generator of a strongly continuous semigroup 
of operators $(T_t)_{t\geq 0}$ given by the formula
$$T_t\varphi=\lim_{n\to\infty}\left(V_{\frac{t}{n}}\right)^n\varphi$$
where the limit exists for every  $\varphi\in \mathcal{X}$ and is uniform  with respect to $t\in [0,t_0]$ for every $t_0>0$. 
Moreover $(T_t)_{t\geq 0}$ satisfies the estimate $\|T_t\|\leq M e^{\omega t}$ for every $t\geq 0$.

\thm\label{TC}(\textit{Approximation of generator implies approximation of semigroup}) (theorem 4.9 in \cite{EN})

Let $(e^{\overline{\mathcal{L}_j}t})_{t\geq 0}$ be a sequence of strongly continuous semigroups of operators in a Banach space $\mathcal{X}$  with the
generators $(\overline{\mathcal{L}_j},Dom(\overline{\mathcal{L}_j}))$, which satisfies, for some fixed constants
$M\geq 1, w\in\R$, the condition $\left\|e^{\overline{\mathcal{L}_j}t}\right\|\leq Me^{wt}$ for all $t\geq 0$ and every $j\in\mathbb{N}$. Suppose there is
a closed linear operator  $(\mathcal{L},Dom(\mathcal{L}))$
on $\mathcal{X}$  with a dense domain $Dom(\mathcal{L})$, such that $\overline{\mathcal{L}_j}x\to Lx$ for every $x\in Dom(\mathcal{L})$.
Suppose the image of the operator $(\lambda_0 I-\mathcal{L})$ is dense in $\mathcal{X}$ for some $\lambda_0>0.$

Then the semigroups $(e^{\overline{\mathcal{L}_j}t})_{t\geq 0}, j\in\N$  converge strongly
(and uniformly in $t\in[0,t_0]$ for every fixed $t_0>0$)
to a strongly continuous semigroup $(e^{\overline{\mathcal{L}}t})_{t\geq 0}$ with the generator $\overline{\mathcal{L}}$.
In other words,
for every $x\in \mathcal{X}$ there exists $\lim_{j\to\infty}e^{\overline{\mathcal{L}_j}t}x=e^{\overline{\mathcal{L}}t}x$ uniformly in $t\in[0,t_0]$ for every fixed $t_0>0$.

\rem Below, the role of  $\mathcal{X}$ will be played by space 
$X$, a closed real subspace of the complex Banach space  $C_b(H,\mathbb{C})$. Because all the operators used in this paper below 
are real, and (as it will be proven further in theorems \ref{preChern2} and \ref{gDeltaA_thm}) $X$ is invariant with respect to them, 
the above theorems about $\mathcal{X}$ are applicable to $X$.

\subsection{Properties of spaces $D$, $X$, $D_1$}

\rem It directly follows from the definitions of these spaces that

i) $D\subset D_1\subset X\subset C_b(H,\R)\subset C_b(H,\mathbb{C})$;

ii) $D$ and $D_1$ are dense in $X$;

iii) $X$ is a Banach space.

\pr\label{prravnneprD} If $f\in D$, then $f$ is uniformly continuous.

\textbf{Proof.} It follows from the definition of the space $D$ that the function $D\ni f\colon H\to\R$ is bounded and 
its Fr\'{e}chet derivatives of all orders exist and are bounded. In particular, there exists

\be\label{ogrpr} \sup_{x\in H}\|f'(x)\| =M<\infty.\ee

One can set $n=0$ in Taylor's formula (\ref{Teyl}) to ensure that for every $x\in H$ and every  $y\in H$ there exists a real number $R_1(x,y)$ such that

\be\label{razlTejl2}f(x)-f(y)=R_1(x,y),\ee

and the estimate holds

\be\label{ocostchl2}|R_1(x,y)| \stackrel{(\ref{oc_ost_chl})}{\leq} \frac{\|y-x\|^{1}}{1!}\sup\limits_{z\in[x,y]}\left\|f'(z)\right\|\stackrel{(\ref{ogrpr})}{\leq} M\|x-y\|.\ee

Hence, for each $x\in H$ and each $y\in H$ we have

\be\label{fravnomnepr}|f(x)-f(y)|\stackrel{(\ref{razlTejl2})}{=}\|R_1(x,y)\|\stackrel{(\ref{ocostchl2})}{\leq} M\|x-y\|,\ee
which implies the uniform continuity of $f$.\QED

\pr\label{prravnneprX} If $\varphi\in X$, then $\varphi$ is uniformly continuous.

\textbf{Proof.} Take any given $\varepsilon>0$. Let us find $\delta>0$ such that $\|x-y\|<\delta$ implies $|\varphi(x)-\varphi(y)|<\varepsilon$.

As $\varphi\in X$, there exists a sequence of functions $(f_j)\subset D$ converging to $\varphi$ uniformly. Hence, there exists a number $j_0$ such
that (introducing the notation $f_{j_0}=f$) we have
\be\label{eps3oc} \left\|\varphi-f_{j_0}\right\|=\left\|\varphi-f\right\|=\sup_{x\in H}|\varphi(x)- f(x)|<\frac{\varepsilon}{3}.\ee

Moreover, as $f\in D$, proposition \ref{prravnneprD} implies estimate (\ref{fravnomnepr}) with some $M>0$.

Let us set $\delta=\frac{\varepsilon}{3M}$ and note that $\|x-y\|<\delta.$ Then
$$|\varphi(x)-\varphi(y)|\leq |\varphi(x)-f(x)|+|f(x)-f(y)|+|f(y)-\varphi(y)| 
\stackrel{(\ref{fravnomnepr}),(\ref{eps3oc})}{<}\frac{\varepsilon}{3} + M\frac{\varepsilon}{3M}+\frac{\varepsilon}{3}=\varepsilon.$$
\QED

\pr\label{ravnomneprsemejstv} Suppose that a sequence of functions  $(f_j)_{j=1}^\infty\subset X$ converges 
uniformly to a function $f_0\in X$. Then the family $(f_j)_{j=0}^\infty$ is equicontinuous.

\textbf{Proof.} Suppose $\varepsilon>0$ is given. Let us find $\delta>0$ such that $\|x-y\|<\delta$ implies that $|\varphi(x)-\varphi(y)|<\varepsilon$.

By proposition \ref{prravnneprX}, function $f_j$ is uniformly continuous for each $j=0,1,2,\dots$. Thus, for each  $j=0,1,2,\dots$ there exists $\delta_j>0$ 
such that $\|x-y\|<\delta_j$ implies
\be\label{deltajoc}|f_j(x)-f_j(y)|<\frac{\varepsilon}{3}.\ee

As $f_j\to f_0$ uniformly, there exists $j_0$ such that for all $j>j_0$

\be\label{ravnomotkl}\sup_{x\in H}|f_0(x)-f_j(x)|<\frac{\varepsilon}{3}.\ee

Let us set  $\delta=\min(\delta_0,\delta_1,\dots,\delta_{j_0}).$ Then for $j>j_0$ we have that $\|x-y\|<\delta$ implies
$$|f_j(x)-f_j(y)|\leq |f_j(x)-f_0(x)|+|f_0(x)-f_0(y)|+|f_0(y)-f_j(y)|\stackrel{(\ref{deltajoc}),(\ref{ravnomotkl})}{<}\frac{\varepsilon}{3}+
\frac{\varepsilon}{3}+\frac{\varepsilon}{3}=\varepsilon.$$
Now, if $0\leq j\leq j_0$, then $\|x-y\|<\delta$ implies estimate (\ref{deltajoc}), which is even stronger.
\QED

\rem A number $a\in \R$ is called a limit at infinity of a funtion $f\colon H\to\R$ if 
$$\lim_{R\to+\infty}\sup_{\|x\|\geq R}|f(x)-a|=0.$$ It is shown in \cite{Remizov} that if $H$ is infinite-dimensional, then 
a non-constant function that belongs to $X$ cannot have a limit at infinity. For example, the function $x\longmapsto \exp(-\|x\|^2)$ 
belongs to $C_b(H,\R)$ but not to $X$.

\rem Suppose that $\alpha_k\colon \R\to\R$ is a family of infinitely-smooth functions, uniformly bounded with their first and second derivatives:
$$\sup_{p\in\{0,1, 2\}}\sup_{k\in\N}\sup_{t\in\R}\left|\frac{d^p\alpha_k(t)}{dt^p}\right|\leq M\equiv\textrm{const}.$$
For example, $\alpha_k(t)=\sin(d_k(t-t_k))$, where $d_k$ and $t_k$ are constants and
$0<d_k\leq 1$. Suppose numerical series $\sum_{k=1}^\infty b_k$ converges absolutely.
Let $(e_k)_{k=1}^\infty$ be an orthonormal basis in $H$.

Then function $$f(x)= \sum_{k=1}^\infty b_k \alpha_k(\left<x,e_k\right>)$$
belongs to the class $D_1$.

This statement can be easily extended to the case $\alpha_k\colon \R^{n_k}\to \R$.

\rem\label{non-separ-D} Space $D$ is not separable (it does not have a countable dense subset). In the case of one-dimensional $H$ it can be 
shown similar to the standard proof of the nonseparability of $C_b(\R,\R)$. 
If $\mathrm{dim} H>1$, then $\R^1$ can be embedded into $H$ as a linear span of 
a non-zero vector $e\in H$. Using this, one can embed the set of cylindrical functions contributing 
to the non-separability of $D$ in the case of one-dimensional $H$, into the space $D$ in the general case.

\rem By Remark \ref{non-separ-D} and the inclusion $D\subset D_1\subset X$, one can see that $D_1$ and $X$ are not separable too.

\section{Main results}

\subsection{Family $S_t$ provides a semigroup with generator $\overline{L}$}

\thm\label{preChern2} (\textit{On the properties of family $(S_t)_{t\geq 0}$
and its connection to the operator $L$})

Suppose that $g\in X$, and for every $x\in H$ we have $g(x)\geq g_o\equiv\mathrm{const}>0$. Suppose 
that $B\in X_H$ and $C\in X$. Suppose that $t>0$, and $\mu_{2tg(x)A}$ is the centered Gaussian 
measure on $H$ with the correlation operator $2tg(x)A$.

For $t\geq 0$ and $\varphi\in C_b(H,\R)$ let us define
\be\label{S_t_def} (S_t\varphi)(x):= e^{tC(x) - t\frac{\langle AB(x),B(x) \rangle}{g(x)}}\int_H\varphi(x+y)e^{\left<\frac{1}{g(x)}B(x),y\right>}\mu_{2tg(x)A}(dy) \textrm{ for }t>0, 
\textrm{ and }S_0\varphi:=\varphi.\ee

Then:

1. If $t\geq 0$ and $\varphi\in C_b(H,\R)$, then $S_t\varphi\in C_b(H,\R)$. For every $t\geq 0$
the operator $S_t\colon C_b(H,\R)\to C_b(H,\R)$ is linear and bounded; its norm does not exceed $e^{\left(\frac{2\|A\|\|B\|^2}{g_0}+\|C\|\right)t }$.

2. If $g\in D$, $C\in D$, $B\in D_H$, then the space $D$ for every $t\geq 0$ is invariant with respect to the operator $S_t$.

3. If $g\in X$, $C\in X$, $B\in X_H$, then the space $X$ for every $t\geq 0$ is invariant with respect to the operator $S_t$.

4. For every function $\varphi\in D$, for $g\in X$, $C\in X$, $B\in X_H$ there exists (uniformly with respect to $x\in H$) a limit

$$\lim\limits_{t\to 0} \frac{(S_t\varphi)(x) -\varphi(x)}{t}=g(x)\tr A\varphi''(x)+\langle \varphi'(x),AB(x)\rangle + C(x)\varphi(x)=(L\varphi)(x).$$

5. If $\varphi\in X$, $g\in X$, $C\in X$, $B\in X_H$, then the function $[0,+\infty) \ni t\longmapsto S_t\varphi\in X$ is 
continuous, i.e. if $t_0\geq 0, t_n\geq 0$ and $t_n\to t_0$, then $\sup_{x\in H}|(S_{t_n}\varphi)(x)-(S_{t_0}\varphi)(x)|\to 0.$

\textbf{Proof.}

1. Function $\varphi$ is bounded, so integral (\ref{S_t_def}) exists by lemma \ref{polyexpint}. Suppose a number $t>0$ and a 
function $\varphi\in C_b(H,\R)$ are fixed. Recalling lemma \ref{linzam}, one can see that
$$\int_H \varphi(x+y)e^{\left<\frac{1}{g(x)}B(x),y\right>}\mu_{2tg(x)A}(dy)=
\int_H\varphi\left(x+\sqrt{2tg(x)}y\right)e^{\left<\frac{1}{g(x)}B(x),\sqrt{2tg(x)}y\right>}\mu_{A}(dy).$$

Introducing the notation $\|B\|=\sup_{x\in H}\|B(x)\|$, we obtain the estimate 
$$\|S_t\varphi\|= \sup_{x\in H}\left|e^{tC(x)- \frac{\langle AB(x),B(x) \rangle}{g(x)}t} \int_H\varphi\left(x+\sqrt{2tg(x)}y\right)e^{\left<\frac{1}{g(x)}B(x),\sqrt{2tg(x)}y\right>}\mu_{A}(dy)\right|\leq$$

$$ \sup_{x\in H}\left|e^{tC(x)- \frac{\langle AB(x),B(x) \rangle}{g(x)}t}\right|\ \ \sup_{x\in H}|\varphi(x)|\ \sup_{x\in H}\int_He^{\left<\sqrt{\frac{2t}{g(x)}}B(x),y\right>}\mu_{A}(dy)\stackrel{(\ref{int_eksp})}{=}$$

\be\label{ocenka_normy_St}e^{t\left(\|C\|+ \left\| \frac{\langle AB(x),B(x) \rangle}{g(x)}\right\|\right)} \|\varphi\| \sup_{x\in H}e^{\frac{1}{2}\frac{2t}{g(x)}\left<AB(x),B(x)\right>}\stackrel{\left<x_1,x_2\right>\leq\|x_1\|\cdot\|x_2\|}{\leq}e^{\left(\frac{2\|A\|\|B\|^2}{g_0}+\|C\|\right)t}\|\varphi\|, \ee

which implies that the function $x\longmapsto (S_t\varphi)(x)$ is bounded. Let us prove that this function is continuous.
Suppose $x_j\to x$, then for every $y\in H$ 
$$\varphi\left(x_j+\sqrt{2tg(x_j)}y\right)e^{\left<\sqrt{\frac{2t}{g(x_j)}}B(x_j),y\right>}\to \varphi\left(x+\sqrt{2tg(x)}y\right)e^{\left<\sqrt{\frac{2t}{g(x)}}B(x),y\right>}.$$

Moreover, 
$\left|\varphi\left(x_j+\sqrt{2tg(x_j)}y\right)e^{\left<\sqrt{\frac{2t}{g(x_j)}}B(x_j),y\right>}\right|\leq \|\varphi\| e^{\sqrt{\frac{2t}{g_0}}\|B\|\|y\|}$
and, in a similar way,

$$\left|\varphi\left(x_j+\sqrt{2tg(x)}y\right)e^{\left<\sqrt{\frac{2t}{g(x)}}B(x),y\right>}\right|\leq \|\varphi\| e^{\sqrt{\frac{2t}{g_0}}\|B\|\|y\|}.$$ Lemma \ref{polyexpint} implies that the function $y\longmapsto e^{\sqrt{\frac{2t}{g_0}}\|B\|\|y\|}$ is integrable over 
the measure $\mu_A$. Therefore by Lebesgue dominated convergence theorem

$$\lim_{j\to\infty}\int_H\varphi\left(x_j+\sqrt{2tg(x_j)}y\right)e^{\left<\sqrt{\frac{2t}{g(x_j)}}B(x_j),y\right>}\mu_{A}=\int_H\varphi\left(x+\sqrt{2tg(x)}y\right)e^{\left<\sqrt{\frac{2t}{g(x)}}B(x),y\right>}\mu_{A}.$$

Because the functions $C, B, g$ are continuous and $g(x)\geq g_0>0$, we have $e^{tC(x_j)- \frac{\langle AB(x_j),B(x_j) \rangle}{g(x_j)}t}\to e^{tC(x)- \frac{\langle AB(x),B(x) \rangle}{g(x)}t}$. Therefore, $(S_t\varphi)(x_j)\to (S_t\varphi)(x)$. 
So we have proved that $S_t\varphi\in C_b(H,\R)$. Estimate (\ref{ocenka_normy_St}) shows that $\|S_t\|\leq e^{\left(\frac{2\|A\|\|B\|^2}{g_0}+\|C\|\right)t }$.\\

2. We fix $t>0$ and prove that $S_t\varphi\in D.$

i) First of all, if $g\in D$, $C\in D$, $B\in D_H$, then the operator $S_t$
maps the cylindrical function $\varphi$ into a cylindrical function $S_t\varphi$. This follows from the fact that (\ref{S_t_def}) 
is a cylindrical function of cylindrical functions, which are functions of finite number of linear functionals of $x$. 
Therefore, the number $(S_t\varphi)(x)$ depends on $x$ only via a finite number of linear functionals, hence $x\longmapsto (S_t\varphi)(x)$ is a cylindrical function, see (\ref{St_cyl-1}) for the exact formula. 

ii) Let us introduce some notation. As $\varphi$ is cylindrical, then for every  $x\in H$
the equality $\varphi(x)=\varphi^n(\langle x,e_1\rangle,\dots,\langle x,e_n\rangle)$ holds
for some $n\in\mathbb{N}$, some function $\varphi^n\colon\mathbb{R}^n\to\mathbb{R}$
and some set of vectors $e_1,\dots,e_n$ of the space $H$. Functions $g$, $C$, $B$ are also cylindrical, 
and without loss of generality we can accept that the set of vectors  $e_1,\dots,e_n$ is so large that the 
following holds:  $g(x)=g^n(\langle x,e_1\rangle,\dots,\langle x,e_n\rangle)$, $C(x)=C^n(\langle x,e_1\rangle,\dots,\langle x,e_n\rangle)$, $B(x)=B_1(\langle x,e_{1}\rangle,\dots,\langle x,e_{n}\rangle)e_1+\dots+B_n(\langle x,e_{1}\rangle,\dots,\langle x,e_{n}\rangle)e_n$.

At this moment vectors $e_1,\dots,e_n$ can be in arbitrary position with respect to the eigenvectors of the operator $A$.
Without loss of generality the set $e_1,\dots,e_n$ can be considered orthonormalized. 

Let us introduce the definitions:  $\Psi_n\colon H\ni h\longmapsto
(\left<h,e_1\right>,\dots,\left<h,e_n\right>)\in\mathbb{R}^n$ --- a projector,
$H_n=\mathrm{span}( e_1,\dots, e_n)$ --- a subspace in $H$, $I_n\colon H_n\ni h\longmapsto
(\left<h,e_1\right>,\dots,\left<h,e_n\right>)\in\mathbb{R}^n$ --- an isomorphism,
$P_n\colon H\ni h\longmapsto \langle h,e_1\rangle e_1+\dots+\langle h,e_n\rangle e_n\in H_n$ --- a projector. Next,
denote $\stackrel{\rightarrow}{x}_1^n=(x_1,\dots,x_n)\in\R^n$ and 
$\stackrel{\rightarrow}{B}_1^n(\stackrel{\rightarrow}{x}_1^n)=(B_1(x_{1},\dots,x_{n}),\dots,B_n(x_{1},\dots,x_{n}))\in \R^n$. With these definitions we have $\Psi_n=I_n P_n$ and $\varphi(x)=\varphi^n(\Psi_n x),$ $g(x)=g^n(\Psi_n x),$ $C(x)=C^n(\Psi_n x),$ $B(x)=\stackrel{\rightarrow}{B}_1^n(\Psi_nx).$

Let us introduce the function $\Phi\colon\R^n\to\R$ in the following way:

\be\label{Phi_def}\Phi\left(\stackrel{\rightarrow}{x}_1^n\right)= \int_H\varphi^n\left(\stackrel{\rightarrow}{x}_1^n+\sqrt{2tg^n\left(\stackrel{\rightarrow}{x}_{1}^n\right)}\Psi_n(y)\right) e^{\sqrt{\frac{2t}{g^n\left(\stackrel{\rightarrow}{x}_1^n \right)}}\left<\stackrel{\rightarrow}{B}_1^n(\stackrel{\rightarrow}{x}_1^n),\Psi_n(y)\right>}  \mu_A(dy).\ee

Then for every $x\in H$ we have 
\be\label{St_cyl-1}(S_t\varphi)(x)=\Phi(\Psi_nx)\exp\left(tC^n(\Psi_nx) -\frac{\left\langle A\stackrel{\rightarrow}{B}_1^n(\Psi_nx),\stackrel{\rightarrow}{B}_1^n(\Psi_nx) \right\rangle}{g^n(\Psi_nx)}t\right).\ee

iii) Now let us prove that $S_t\varphi$ has bounded Fr\'{e}chet derivatives of all orders employing the proposition \ref{pr_diff}. To do this we need to prove that the functions $\R^n\to\R$ have bounded Fr\'{e}chet derivatives of all orders. The exponent in (\ref{St_cyl-1}) has this property because the exponent is obtained by composition and arithmetical operations from the functions with this property.

Let us show that $\Phi$ has Fr\'{e}chet derivatives of all orders. The product of differentiable functions under the sign of integral in (\ref{Phi_def}) is differentiable, so the problem is reduced to the verification of the differentiability of the integral. To do this, we apply 
lemma \ref{integrir_cil_f} and arrive from an integral over $H$
to an integral over $\R^n$ in the expression for $\Phi$ (this is possible because the intergand is cylindrical).

Operator $A$ is non-degenerate and symmetric on $H$, therefore the operator $P_n A$ is non-degenerate and symmetric on $H_n$,
and therefore it can be diagonalized  in some orthonormal basis $b_1,\dots,b_n$. Without loss of generality we can assume that the vectors $e_1,\dots,e_n$ form such a basis. Indeed, changing the basis in the space $H_n$ will just produce linear non-degenerate change of variables in the functions $\R^n\to\R$ used to define cylindrical functions $H\to\R$. This will give us new functions $\R^n\to\R$, but all their properties that we need will be preserved.

The matrix of the operator $P_n A$ in $H_n$ coincides with the matrix of the operator $Q_n=I_nP_nAI_n^{-1}$ in  $\R^n$. Next, let $q_1,\dots,q_n$ be
the eigenvalues of the operator $Q_n$, corresponding to the eigenvectors $\Psi_ne_1,\dots,\Psi_ne_n$. Note, that $q_i>0$ and
$g^n\left(\stackrel{\rightarrow}{x}_{1}^n\right)\geq g_0\equiv\mathrm{const}>0$
for every $\stackrel{\rightarrow}{x}_{1}^n\in\R^n$. Then by (\ref{intcyl}) we have

$$\Phi\left(\stackrel{\rightarrow}{x}_1^n\right)=\left(\frac{1}{\sqrt{2\pi}}\right)^n\frac{1}{\sqrt{\prod_{i=1}^n q_i}} \int_{\mathbb{R}^n} \varphi^n\left(\stackrel{\rightarrow}{x}_1^n+\sqrt{2tg^n\left(\stackrel{\rightarrow}{x}_{1}^n\right)}\ z \right) \times$$

$$e^{\sqrt{\frac{2t}{g^n\left(\stackrel{\rightarrow}{x}_1^n \right)}}\left<\stackrel{\rightarrow}{B}_1^n(\stackrel{\rightarrow}{x}_1^n),z\right>}\exp\left(-\sum_{i=1}^n \frac{z_i^2}{2q_i} \right)dz.$$

Now we introduce a measure $\nu$ on $\R^n$ given by its density with respect to the Lebesgue measure:
for every measurable set $\mathcal{A}\subset\R^n$ we set
$$\nu(\mathcal{A}):=\left(\frac{1}{\sqrt{2\pi}}\right)^n\frac{1}{\sqrt{\prod_{i=1}^n q_i}} \int_{\mathcal{A}}\exp\left(-\sum_{i=1}^n \frac{z_i^2}{2q_i} \right)dz.$$
It follows from the definitions given above that

\be\label{nu_int}\Phi\left(\stackrel{\rightarrow}{x}_1^n\right)=\int_{\mathbb{R}^n} \varphi^n\left(\stackrel{\rightarrow}{x}_1^n+\sqrt{2tg^n\left(\stackrel{\rightarrow}{x}_{1}^n\right)}\ z \right) e^{\sqrt{\frac{2t}{g^n\left(\stackrel{\rightarrow}{x}_1^n \right)}}\left<\stackrel{\rightarrow}{B}_1^n(\stackrel{\rightarrow}{x}_1^n),z\right>}  \nu(dz).\ee

The integrand in (\ref{nu_int}) is a composition of mappings with the continuous bounded Fr\'{e}chet derivative. Thus, it
has a continuous bounded Fr\'{e}chet derivative. The Fr\'{e}chet derivative of the integrand is uniformly bounded (the estimate is obtained from the chain rule formula), and $(\R^n,\nu)$ is locally compact, countable at infinity, linear normed space with the non-negative Radon measure. Therefore we can apply theorem 115 from  \cite{Schwartz}
on the Fr\'{e}chet differentiation  under the Lebesgue integral.
Repeating this reasoning for every $k\in \N$, we conclude that as the integrand has, everywhere in  $\R^{n}$, continuous Fr\'{e}chet derivatives of $k$-th order, 
then the function $\Phi$ has, everywhere in $\R^{n}$,
continuous bounded Fr\'{e}chet's derivatives of $k$-th order. So the functions $\R^n\to\R$ in the right-hand side of (\ref{St_cyl-1}) all have continuous bounded Fr\'{e}chet's derivatives of $k$-th order. 

Therefore, according to point 2 of proposition \ref{pr_diff}, the function $x\longmapsto (S_t\varphi)(x)$ also has, for every $k\in\N$, Fr\'{e}chet derivatives of order $k$, continuous and bounded everywhere on $H$. Therefore $S_t\varphi\in D$.

3. i) Now suppose $\varphi\in X$, which means that $\varphi\in C_{b}(H,\R)$ and there exists a sequence  $(\varphi_j)\subset D$ such that
$\varphi_j(x)\to\varphi(x)$ uniformly with respect to $x\in H$. Suppose also that $g\in X$, so $g\in C_{b}(H,\R)$ and there exists a sequence $(g_j)\subset D$ such
that $g_j\to g$ uniformly. It follows from $g(x)\geq g_0\equiv\mathrm{const}>0$ for all $x\in H$ that there exists a number $j_0\in\N$ such
that for all $j>j_0$ and for all $x\in H$ the inequality $g_j(x)\geq \frac{g_0}{2}$ holds. Therefore, we 
will not restrict the generality when assuming that for the sequence $(g_j)$ the inequality $g_j(x)\geq \frac{g_0}{2}$ already holds for all $j\in\N$ and for all $x\in H$.

Also suppose $C\in X$, so $C\in C_{b}(H,\R)$ and there exists a sequence $(C_j)\subset D$ such that $C_j\to C$ uniformly. Finally, 
suppose $B\in X_H$, so $B\in C_{b}(H,H)$ and there exists a sequence $(B_j)\subset D_H$ such that $B_j\to B$ uniformly. Let $t>0$ be fixed as before.

Let us denote the operator $S_t$ constructed with functions $g_j$, $B_j, C_j$ by $(S_j)_t$. According to (just proven above) item 2 of the theorem,
$(S_j)_t\varphi_j\in D$ for all $j\in\N$. Next, in ii) and iii) we will prove that $((S_{j})_t\varphi_{j})(x)\to (S_t\varphi)(x)$
uniformly with respect to $x\in H$; by the definition of the space $X$ this will mean that $S_t\varphi\in X$.

ii) First of all let us prove that for every fixed $y\in H$ the sequence of functions $x\longmapsto \varphi_j(x+\sqrt{2tg_j(x)}y)$ converges to 
the function $x\longmapsto \varphi(x+\sqrt{2tg(x)}y)$ uniformly with respect to $x\in H$. Indeed, suppose $\varepsilon>0$ is given. 
Let us find $j^*\in\N$ such that for all $j>j^*$ the following estimate holds
\be\label{willthis}\sup_{x\in H}\left|\varphi_j\left(x+\sqrt{2tg_j(x)}y\right)
-\varphi\left(x+\sqrt{2tg(x)}y\right)\right|\leq\varepsilon.\ee

Notice that $\varphi,\varphi_j,g,g_j$ are elements of $X$ by the hypothesis made above, and all the functions in $X$ are uniformly continuous,
according to proposition \ref{prravnneprX}.  According to proposition \ref{ravnomneprsemejstv}, the family of functions $\{\varphi_j:j\in\N\}$ is equicontinuous. 
Thus, there exists $\delta>0$ such that, for all $j\in\N$, if $\|x_1-x_2\|<\delta$, then
\be\label{fjravnostnepr}|\varphi_j(x_1)-\varphi_j(x_2)|<\frac{\varepsilon}{2}.\ee

Function $[0,+\infty)\ni a\longmapsto\sqrt{2ta}\in\R$ is uniformly continuous, so the uniform (with respect to $x\in H$) convergence  $g_j(x)\to g(x)$ implies 
the uniform (with respect to $x\in H$) convergence $[H\ni x\longmapsto \sqrt{2tg_j(x)}\in\R]\to [H\ni x\longmapsto \sqrt{2tg(x)}\in\R]$. This shows
that for every fixed $y\in H$ there exists $j_1\in\N$ such that for all $j>j_1$ and all $x\in H$ we have
\be\label{ocdeltaj1}\left\|\left(x+\sqrt{2tg_j(x)}y\right)-\left(x+\sqrt{2tg(x)}y\right)\right\|<\delta.\ee

Besides, because $\varphi_j(z)\to \varphi(z)$ uniformly with respect to $z\in H$, there exists a number $j_2$ such that for all $j>j_2$ and all $z\in H$ 
we have

\be\label{ravnomshodt}|\varphi_j(z)-\varphi(z)|<\frac{\varepsilon}{2}.\ee

For each fixed $y\in H$, for all $j\in\N$ and for all $x\in H$, we have
$$\left|\varphi_j\left(x+\sqrt{2tg_j(x)}y\right) - \varphi\left(x+\sqrt{2tg(x)}y\right)\right|\leq$$
$$\left|\varphi_j\left(x+\sqrt{2tg_j(x)}y\right) - \varphi_j\left(x+\sqrt{2tg(x)}y\right)\right| + \left|\varphi_j\left(x+\sqrt{2tg(x)}y\right) - \varphi\left(x+\sqrt{2tg(x)}y\right)\right|.$$

Now, we define  $j^*=\max(j_1,j_2).$ For all $j>j^*$ the first summand is less than $\frac{\varepsilon}{2}$ due to (\ref{fjravnostnepr}) and (\ref{ocdeltaj1}); 
the easiest way to see that is to set  $x_1=x+\sqrt{2tg_j(x)}y$ and $x_2=x+\sqrt{2tg(x)}y$. As for the second summand, it is less than 
$\frac{\varepsilon}{2}$ on account of (\ref{ravnomshodt}); the easiest way to see that is to set $z=x+\sqrt{2tg(x)}y$.

So, for each fixed $y\in H$ we have found a number $j^*\in\N$ such that for all $j>j^*$ and for all $x\in H$ the following inequality holds:
$$\left|\varphi_j\left(x+\sqrt{2tg_j(x)}y\right) - \varphi\left(x+\sqrt{2tg(x)}y\right)\right|<\frac{\varepsilon}{2}+\frac{\varepsilon}{2}=\varepsilon.$$

By taking $\sup_{x\in H}$ we obtain the needed estimate (\ref{willthis}), as the right hand side of the inequality above does not depend on $x$.

iii) A reasoning similar to ii) shows that for fixed $y$ we have
$e^{\left<\sqrt{\frac{2t}{g_j(x)}}B_j(x),y\right>}\to e^{\left<\sqrt{\frac{2t}{g(x)}}B(x),y\right>}$ uniformly with respect to $x\in H$; let us omit the 
detailed proof. Uniting this with the results of ii) and keeping in mind that for fixed $y$ all the sequences of functions
are bounded collectively, we obtain that for each fixed  $y$ we have

\be\label{eqpred}\varphi_j\left(x+\sqrt{2tg_j(x)}y\right)e^{\left<\sqrt{\frac{2t}{g_j(x)}}B_j(x),y\right>}\to \varphi\left(x+\sqrt{2tg(x)}y\right)e^{\left<\sqrt{\frac{2t}{g(x)}}B(x),y\right>}\ee

uniformly with respect to $x\in H.$

iv) For fixed $y$, functions in (\ref{eqpred}) are bounded, therefore the sequence of functions $Y_j\colon H\to\R$
$$Y_j=\left[y\longmapsto \sup_{x\in H}\left|\varphi_j\left(x+\sqrt{2tg_j(x)}y\right)e^{\left<\sqrt{\frac{2t}{g_j(x)}}B_j(x),y\right>} -  \varphi\left(x+\sqrt{2tg(x)}y\right)e^{\left<\sqrt{\frac{2t}{g(x)}}B(x),y\right>} \right| \right]$$ is well defined.
It follows from iii), that $Y_j(y)$ converges to zero pointwise, in other words for each $y$. Functions $Y_j$ are non-negative and bounded collectively 
by an integrable (due to lemma \ref{polyexpint}) function $y\longmapsto \alpha e^{\beta\|y\|}+\gamma$ with appropriate constants $\alpha, \beta$ and $\gamma$.
The Lebesgue dominated convergence theorem gives us that $\int_HY_j(y)\mu_A(dy)\to 0$. As the number sequence $\int_HY_j(y)\mu_A(dy)$ converges, it is bounded. 
For brevity, let us denote $$\Psi(x,y)=\varphi\left(x+\sqrt{2tg(x)}y\right)e^{\left<\sqrt{\frac{2t}{g(x)}}B(x),y\right>}, \Psi_j(x,y)=\varphi_j\left(x+\sqrt{2tg_j(x)}y\right)e^{\left<\sqrt{\frac{2t}{g_j(x)}}B_j(x),y\right>},$$

$$E(x)=\exp\left(tC(x) - t\frac{\langle AB(x),B(x) \rangle}{g(x)}\right), E_j(x)=\exp\left(tC_j(x) - t\frac{\langle AB_j(x),B_j(x) \rangle}{g_j(x)}\right).$$

A reasoning similar to ii) shows that $E_j(x)\to E(x)$ uniformly with respect to $x\in H$. We obtain the estimate
$$\left\|(S_j)_t\varphi_j - S_t\varphi\right\|=\sup_{x\in H}\left| E_j(x)\int_H \Psi_j(x,y) \mu_A(dy) -E(x)\int_H \Psi(x,y) \mu_A(dy) \right|\leq$$

$$\sup_{x\in H}\left| E_j(x)\int_H \Psi_j(x,y) \mu_A(dy) -E_j(x)\int_H \Psi(x,y) \mu_A(dy) \right|+$$

$$\sup_{x\in H}\left| E_j(x)\int_H \Psi(x,y) \mu_A(dy) -E(x)\int_H \Psi(x,y) \mu_A(dy) \right|\leq$$

$$\sup_{x\in H} \left| E_j(x)\right| \int_H Y_j(y)\mu_A(dy) + \sup_{x\in H} \left| E_j(x)-E(x)\right| \int_H \Psi(x,y) \mu_A(dy).$$

Finally let us note that $\sup_{x\in H} \left| E_j(x)-E(x)\right|\to 0$ and $\int_H Y_j(y)\mu_A(dy)\to 0$, and the multipliers of these terms in the above 
estimate are bounded, which implies $\left\|(S_j)_t\varphi_j - S_t\varphi\right\|\to 0$.\\

4. Suppose $\varphi\in D$, $t>0$ and $x\in H$ are fixed. Let us consider the integral $$\int_H\varphi(x+y)e^{\left<\frac{1}{g(x)}B(x),y\right>}\mu_{2tg(x)A}(dy).$$

We will work with the approximation of the function $\varphi$ by its Taylor polynomial (\ref{Teyl}) of the second order with the center at the point $x$. 
Before we start, it is important to note that the remainder term $R(x,y)$ will not be small, as the vector $y$ ranges over the whole space $H$,
and vector $x$ is fixed. However as $\varphi$ is three times Fr\'{e}chet differentiable on $H$,
for each $x\in H$ and each $y\in H$ there exists a real number $R(x,y)$ such that

$$\int_H\varphi(x+y)e^{\left<\frac{1}{g(x)}B(x),y\right>}\mu_{2tg(x)A}(dy)=\int_H\Big\{\varphi(x) + \left<\varphi'(x),y\right> + \frac{1}{2!}\left<\varphi''(x)y,y\right>+$$

$$+ R(x,y)\Big\}e^{\left<\frac{1}{g(x)}B(x),y\right>} \mu_{2tg(x)A}(dy),$$

and, according to (\ref{oc_ost_chl}), the estimate

\be\label{ooc}|R(x,y)|\leq\frac{\|y\|^{3}}{(3)!}\sup\limits_{z\in[x,x-y]}\left\|\varphi^{(3)}(z)\right\|\leq \frac{1}{3!} \|\varphi'''\|\|y\|^3\ee

holds, where we define $\|\varphi'''\|=\sup_{z\in H}\|\varphi^{(3)}(z)\|.$ Also, let us denote $\|g\|=\sup_{z\in H}|g(z)|.$

The sum can be integrated termwise as every term can be bounded by a polynomial of $\|y\|$, multiplied by an exponent of $\|y\|$
and such functions are integrable due to lemma \ref{polyexpint}. Let us calculate integrals in the sum and for fixed $\varphi$ 
evaluate the decay rate as $t\to 0$ uniformly with respect to $x\in H$.

$$\int_H \varphi(x)e^{\left<\frac{1}{g(x)}B(x),y\right>}\mu_{2tg(x)A}(dy)=\varphi(x)\int_H e^{\left<\frac{1}{g(x)}B(x),y\right>}\mu_{2tg(x)A}(dy)=$$

 (set $z=\frac{B(x)}{g(x)}$ and $\widetilde{A}=2tg(x)A$  in (\ref{int_eksp}))

$$\varphi(x)\exp\left(\frac{1}{2}\left\langle 2tg(x)A\frac{B(x)}{g(x)},\frac{B(x)}{g(x)} \right\rangle \right)=\varphi(x)\exp\left(\frac{\left\langle AB(x),B(x) \right\rangle}{g(x)}t \right)=$$

$$\varphi(x)\left(1+ \frac{\langle AB(x),B(x) \rangle}{g(x)}t +o(t) \right)= \varphi(x) + \varphi(x)\frac{\langle AB(x),B(x) \rangle}{g(x)}t +o(t). $$

Next, let $z$ and $\widetilde{A}$ be as before, and set $w=\varphi'(x)$ in (\ref{int_eksp1}). 

$$\int_H \langle \varphi'(x),y\rangle e^{\left<\frac{1}{g(x)}B(x),y\right>}\mu_{2tg(x)A}(dy)\stackrel{(\ref{int_eksp1})}{=}$$

$$\left\langle 2tg(x)A\varphi'(x),\frac{1}{g(x)}B(x)\right\rangle \exp\left(\frac{1}{2}\left\langle 2tg(x)A\frac{1}{g(x)}B(x),\frac{1}{g(x)}B(x) \right\rangle \right)=$$

$$2t\langle A\varphi'(x),B(x)\rangle \exp\left(\frac{t}{g(x)}\left\|\sqrt{A}B(x)\right\|^2\right)=2t\langle A\varphi'(x),B(x)\rangle + o(t).$$

For the term with $\varphi''$ we have

$$\int_H \left<\varphi''(x)y,y\right>e^{\left<\frac{1}{g(x)}B(x),y\right>}\mu_{2tg(x)A}(dy)\stackrel{(\ref{int_eksp2})}{=}$$

$$=\Big(2tg(x)\tr A\varphi''(x) + 4t^2\langle \varphi''(x)AB(x),AB(x)\rangle\Big)\exp\left(\frac{\left\langle AB(x),B(x) \right\rangle}{g(x)}t \right)=2tg(x)\tr A\varphi''(x)+o(t).$$

Finally, $$\left|\int_H R(x,y) e^{\left<\frac{1}{g(x)}B(x),y\right>}\mu_{2tg(x)A}(dy)\right| \stackrel{(\ref{lz})}{\leq} \int_H \left|R(x,\sqrt{2tg(x)}y)\right|e^{\left<\frac{1}{g(x)}B(x),\sqrt{2tg(x)}y\right>} \mu_{A}(dy)\stackrel{(\ref{ooc})}{\leq}$$

$$ \int_H \frac{1}{3!}\|\varphi'''\| \left(\sqrt{2tg(x)}\right)^3\|y\|^3 e^{\left<\sqrt{\frac{2t}{g(x)}}B(x),y\right>} \mu_A(dy)\leq$$

$$t^{\frac{3}{2}}\left(\frac{\sqrt{2}}{3}\|\varphi'''\| \|g\|^{\frac{3}{2}}\int_H\|y\|^3 e^{\left<\sqrt{\frac{2t}{g(x)}}B(x),y\right>} \mu_A(dy)\right)\leq$$

(due to the inequality $\langle z,y\rangle\leq \|z\|\|y\|$ and monotonicity of the exponent function)

$$t^{\frac{3}{2}}\left(\frac{\sqrt{2}}{3}\|\varphi'''\| \|g\|^{\frac{3}{2}}\int_H\|y\|^3 e^{\left\|\sqrt{\frac{2t}{g(x)}}B(x)\right\| \left\|y\right\|} \mu_A(dy)\right)\leq$$

(due to $\|B(x)\|\leq B_0$, $g_0\leq g(x)$ and monotonicity of the exponent function)

$$t^{\frac{3}{2}}\left(\frac{\sqrt{2}}{3}\|\varphi'''\| \|g\|^{\frac{3}{2}}\int_H\|y\|^3 e^{\sqrt{\frac{2t}{g_0}}B_0 \left\|y\right\|} \mu_A(dy)\right)=o(t)\textrm{ uniformly with respect to }x\in H.$$

Summing everything up, one can see that 

$$\int_H\varphi(x+y)e^{\left<\frac{1}{g(x)}B(x),y\right>}\mu_{2tg(x)A}(dy)=\varphi(x) + t\langle \varphi'(x),AB(x)\rangle+tg(x)\tr A\varphi''(x)- t\varphi(x)\frac{\langle AB(x),B(x) \rangle}{g(x)}+o(t),$$
uniformly with respect to $x\in H$.

Consider the term $\exp\left(tC(x) - \frac{\langle AB(x),B(x) \rangle}{g(x)}t\right)$. As $\|C\|, \|A\|, \|B\|$ are bounded from infinity and $g$ is bounded from zero, we have $e^{tC(x) - \frac{\langle AB(x),B(x) \rangle}{g(x)}t}=1+ tC(x) - \frac{\langle AB(x),B(x) \rangle}{g(x)}t +o(t)$ uniformly with respect to $x\in H$. Therefore, 
uniformly with respect to $x\in H$ for $t\to 0$ one obtains

$$(S_t\varphi)(x)=e^{tC(x) - \frac{\langle AB(x),B(x) \rangle}{g(x)}t}\int_H\varphi(x+y)e^{\left<\frac{1}{g(x)}B(x),y\right>}\mu_{2tg(x)A}(dy)=$$

$$\varphi(x) + t\langle \varphi'(x),AB(x)\rangle+tg(x)\tr A\varphi''(x)+ tC(x)\varphi(x)+ o(t).$$

This implies that uniformly with respect to $x\in H$

$$\lim\limits_{t\to 0} \frac{(S_t\varphi)(x) -\varphi(x)}{t}=g(x)\tr A\varphi''(x)+\langle \varphi'(x),AB(x)\rangle + C(x)\varphi(x)=(L\varphi)(x).$$

5. i) First let us consider the case $t_0\neq 0.$ If $t_n\to t_0$, then $2t_ng(x)\to 2t_0g(x)$ uniformly with respect to $x\in H$. Because the function $z\longmapsto \sqrt{z}$ is uniformly continuous, 
it follows that
$\sqrt{2t_ng(x)}\to \sqrt{2t_0g(x)}$ uniformly with respect to $x\in H$.

According to proposition \ref{prravnneprX}, the function $\varphi$ is uniformly continuous. Therefore, for every fixed $y\in H$ 
$\varphi(x+\sqrt{2t_ng(x)}y)\to \varphi(x+\sqrt{2t_0g(x)}y)$ uniformly with respect to $x\in H$.

Next, for every $y\in H$ the sequence $\left<\frac{1}{g(x)}B(x),\sqrt{2t_ng(x)}y\right>$ converges to $\left<\frac{1}{g(x)}B(x),\sqrt{2t_0g(x)}y\right>$ 
uniformly with respect to $x\in H$ because of the linearity of the scalar product. Since the number sequence $t_n$ converges, it is bounded; besides, 
$g(x)\geq g_0\equiv\mathrm{const}>0$ and functions $x\longmapsto g(x)$ and $x\longmapsto \|B(x)\|$ are bounded. Therefore, there exists a constant $K>0$ such
that for fixed  $y\in H$ for all $k=0,1,2,3,\dots$ and all $x\in H$ the inequality 
$\left|\left<\frac{1}{g(x)}B(x),\sqrt{2t_kg(x)}y\right>\right|\leq K$ holds. Function $z\longmapsto e^z$  is uniformly continuous with respect 
to $z\in [-K,K]$, therefore for every $y\in H$ we have the convergence 
$e^{\left<\frac{1}{g(x)}B(x),\sqrt{2t_ng(x)}y\right>}\to e^{\left<\frac{1}{g(x)}B(x),\sqrt{2t_0g(x)}y\right>},$
uniformly with respect to $x\in H$.

The product of two collectively bounded uniformly converging sequences is a sequence, uniformly converging to the product of the limits of these sequences. 
Therefore, from the two last paragraphs, it follows that for every $y\in H$ we have  
\be\label{eqpred1}\varphi(x+\sqrt{2t_ng(x)}y)e^{\left<\frac{1}{g(x)}B(x),\sqrt{2t_ng(x)}y\right>}\to   \varphi(x+\sqrt{2t_0g(x)}y)e^{\left<\frac{1}{g(x)}B(x),\sqrt{2t_0g(x)}y\right>}\ee uniformly with respect to $x\in H$.

Since for fixed $y$ the functions from (\ref{eqpred1}) are bounded, the sequence of the number functions is well defined
$$Y_n=\left[y\longmapsto \sup_{x\in H}\left|\varphi\left(x+\sqrt{2t_ng(x)}y\right)e^{\left<\sqrt{\frac{2t_n}{g(x)}}B(x),y\right>} -  \varphi\left(x+\sqrt{2t_0g(x)}y\right)e^{\left<\sqrt{\frac{2t_0}{g(x)}}B(x),y\right>} \right| \right].$$
According to above results, $Y_n(y)$ converges to zero pointwise, i.e. for every $y$. Functions $Y_n$ are bounded by an integrable function (see lemma \ref{polyexpint}),
therefore, by the Lebesgue dominated convergence theorem, we have  $\int_HY_n(y)\mu_A(dy)\to 0$.

Finally, $C\in X, B\in X_H$, and therefore $\sup_{x\in H}|C(x)|<\infty$, $\sup_{x\in H}\|B(x)\|<\infty$. United with $g(x)\geq g_0$ this implies $e^{t_nC(x)-\frac{\langle AB(x),B(x) \rangle}{g(x)}t_n}\to e^{t_0C(x)- \frac{\langle AB(x),B(x) \rangle}{g(x)}t_0}$ uniformly with respect to $x\in H$.

So, the sequence $S_{t_n}\varphi(x)$ is a product of two collectively bounded uniformly converging sequences; 
thus it converges to the product of the limits of these sequences, i.e. to $S_{t_0}\varphi(x)$.

ii) Now let us consider the case $t_0=0.$ Recall that $S_0=I,$ so we need for each fixed $\varphi$ prove that if $t_n\to 0$ then $S_{t_n}\varphi \to \varphi,$ i.e. $\|S_{t_n}\varphi - \varphi\|\to 0.$ Without loss of generality we can assume that $0<t_n\leq 1.$

iia) Let us first consider the case $\varphi\in D.$ Then by item 4 of the present theorem $S_t\varphi=\varphi +tL\varphi +o(t),$ so $t_n\to 0$ implies $S_{t_n}\varphi\to\varphi$.

iib) Now suppose that $\varphi\in X,$ so there exists a sequence $(\varphi_k)\subset D$ such that $\|\varphi_k-\varphi\|\to 0.$ Item 1 of the present theorem implies that there exists such a constant $\omega\geq 1$ such that for all $t\in[0,1]$ we have $\|S_t\|\leq \omega.$ Now for arbitrary $\varepsilon>0$ we apply "$\varepsilon/3$-method" based on the estimate
$$\|S_{t_n}\varphi - \varphi\|\leq \|S_{t_n}\varphi - S_{t_n}\varphi_k\| + \|S_{t_n}\varphi_k - \varphi_k\| + \|\varphi_k - \varphi\|.$$
We can select such $k$ that $\|\varphi_k-\varphi\|<\varepsilon/(3\omega).$ For fixed $k$ we can thanks to iia) find such a number $n_0$ that for all $n>n_0$ one has $\|S_{t_n}\varphi_k - \varphi_k\|<\varepsilon/3.$ So for all $n>n_0$ we see that
$$\|S_{t_n}\varphi - \varphi\|\leq \|S_{t_n}\|\cdot\|\varphi - \varphi_k\| + \|S_{t_n}\varphi_k - \varphi_k\| + \|\varphi_k - \varphi\|< \omega \frac{\varepsilon}{3\omega} + \frac{\varepsilon}{3}+\frac{\varepsilon}{3\omega}\leq\varepsilon.$$

  \QED

\thm\label{gDeltaA_thm} (\textit{On the properties of the operator $L$})
Suppose for each $x\in H$ the inequalities $g(x)\geq g_o\equiv\mathrm{const}>0$ and $C(x)\leq 0$ hold. As $C\in X$, there exists a sequence 
$(C_j)\subset D$ converging to $C$ uniformly; let us additionally require that this sequence can be selected in such a way that $C_j(x)\leq 0$ 
for all $j\in\N$ and all $x\in H$. The operator $L\colon D\to X$ is defined by the equation
$$(L\varphi)(x)=g(x)\tr A\varphi''(x)+\langle \varphi'(x),AB(x)\rangle + C(x)\varphi(x).$$ Symbol $I$ stands for the identity operator.

Then:

1. If $g\in D$, $C\in D$, $B\in D_H$ and $\varphi\in D$, then $L\varphi\in D.$ If $g\in X$, $C\in X$, $B\in X_H$ and $\varphi\in D$, then $L\varphi\in X.$

2. If $g\in D$, $C\in D$, $B\in D_H$, then for each $\lambda>0$ the operator $\lambda I-L$ is surjective on $D,$ therefore $(\lambda I-L)(D)=D$ is a dense subspace in $X$.

3. If $g\in D$, $C\in D$, $B\in D_H$, then the operator $(L,D)$ is dissipative and closable.

4. If $g\in X$, $C\in X$, $B=0$, then for each  $\lambda>0$ the space  $(\lambda I-L)(D)$ is dense in $X$.

5. If $g\in X$, $C\in X$, $B\in X_H$, then the operator $(L,D)$ is dissipative and has the closure $(\overline{L}, D_1)$. This operator is also dissipative.\\

\textbf{Proof.}

1. First part of the statement is obvious: the result of the applying a differential operator with 
smooth cylindrical coefficients to a smooth cylindrical function is a smooth cylindrical function. As it follows from the chain rule, 
all the derivatives are bounded. Let the coefficients of the operator $L$ be uniform limits  $g, B, C$ of cylindrical functions $g_j, B_j, C_j$. 
We denote the operator $L$ that corresponds to the functions $g_j, B_j, C_j$ as $L_j$. Then, as $j\to\infty$ the sequence $L_j\varphi$ converges uniformly to 
$L\varphi$, therefore $L\varphi\in X$, because $L_j\varphi\in D$ as we proved above.

2. Suppose $g\in D$, $B\in D_H$, $C\in D$. Recall that all these functions are cylindrical and thus they are closely related to functions on $\R^n$; 
this will be the core idea of the proof that follows.

Let us fix $\lambda>0$, choose arbitrary function $\varphi\in D$ and then show, that there exists a function $f\in D$ satisfying the equation
\be\label{plotn} \lambda f(x)-g(x)\tr Af''(x) - \langle f'(x),AB(x)\rangle - C(x)f(x)=\varphi(x).\ee

Let the vectors $e_1,\dots,e_n$ be such that for every $x\in H$ we have

\be\label{cilfun12}C(x)=C^n(\langle x,e_1\rangle,\dots,\langle x,e_n\rangle), B(x)=\sum_{k=1}^n B^n_k(\langle x,e_1\rangle,\dots,\langle x,e_n\rangle)e_k,\ee

\be\label{cilfun1}\varphi(x)=\varphi^n(\langle x,e_1\rangle,\dots,\langle x,e_n\rangle) \textrm{ and }  g(x)=g^n(\langle x,e_1\rangle,\dots,\langle x,e_n\rangle),\ee

where $\varphi^n\colon \R^n\to \R$, $g^n\colon \R^n\to \R$, $C^n\colon \R^n\to \R$ and $B^n_k\colon \R^n\to \R$ are continuously differentiable functions, 
bounded along with all the derivatives. Let us find a solution of the equation (\ref{plotn}) in the form

\be\label{cilfun11}f(x)=f^n(\langle x,e_1\rangle,\dots,\langle x,e_n\rangle),\ee

where $f^n\colon \R^n\to \R$ is a continuously differentiable function, bounded along with all its derivatives.

According to item 1 of proposition \ref{pr_diff}

$$f'(x)=\sum_{s=1}^n\partial_sf^n(\langle x,e_1\rangle,\dots,\langle x,e_n\rangle)e_s.$$

Using (\ref{cilfun12}) we obtain that

$$\langle f'(x),AB(x)\rangle=\left< \sum_{s=1}^n\partial_sf^n(\langle x,e_1\rangle,\dots,\langle x,e_n\rangle)e_s, A \sum_{k=1}^n B^n_k(\langle x,e_1\rangle,\dots,\langle x,e_n\rangle)e_k\right>= $$

\be\label{lab12}\sum_{s=1}^n \left(\sum_{k=1}^n B^n_k(\langle x,e_1\rangle,\dots,\langle x,e_n\rangle)\left<Ae_s,  e_k\right>\right)   \partial_sf^n(\langle x,e_1\rangle,\dots,\langle x,e_n\rangle).\ee

One can see, by taking into account
(\ref{traf}), (\ref{cilfun1}), (\ref{cilfun12}), (\ref{cilfun11}) and (\ref{lab12}), that equation  (\ref{plotn}) for the unknown function $f$ 
of form (\ref{cilfun11}) is equivalent to the equation for an unknown function $f^n$:

$$g^n(x_1,\dots,x_n)\sum_{s=1}^n\sum_{k=1}^n \langle Ae_s,e_k\rangle \partial_k\partial_s f^n(x_1,\dots,x_n)$$

$$+\Bigg(\sum_{s=1}^n \left(\sum_{k=1}^n B^n_k(x_1,\dots,x_n)\left<Ae_s,  e_k\right>\right) \partial_sf^n(x_1,\dots, x_n)\Bigg)$$

\be\label{plotn3} + C^n(x_1,\dots,x_n)f^n(x_1,\dots, x_n) - \lambda f^n(x_1,\dots,x_n) =-\varphi^n(x_1,\dots,x_n),\ee

where we denote $x_j=\langle x,e_j\rangle$.

The equation  (\ref{plotn3}) is a finite-dimensional PDE. Note that $g_n(x_1,\dots,x_n)\geq g_0\equiv\textrm{const}>0$.
The quadratic form given by the $n\times n$-matrix $(\langle Ae_s,e_k\rangle)$ is positively defined because the operator $A$ is positive and non-degenerate. Therefore,
the differential operator in equation (\ref{plotn3}) is elliptic, and (\ref{plotn3}) has form (\ref{urchp}), i.e.

$$L^nf^n(x_1,\dots,x_n) - \lambda f^n(x_1,\dots,x_n)=-\varphi(x_1,\dots,x_n).$$

Functions $g^n$, $C^n$ and $B^n_k$ used to construct the coefficiences of the operator $L^n$, and the function $\varphi^n$  are bounded along with all the derivatives.
(Note that $L^n$ is not the $n$-th power of the operator $L$.)
Therefore we can apply item 1 of lemma \ref{lm_urchp} to equation (\ref{plotn3}), so (\ref{plotn3}) has the solution $f^n$, continuous and bounded along
with all the derivatives. Function $f$ defined by equation (\ref{cilfun11}) belongs to $D$ according to item 2 of proposition \ref{pr_diff}.

Thus, for every $\lambda>0$ the operator $\lambda I-L$ is surjective in $D$,
because in the preimage of the function  $\varphi\in D$ there is at least the function
$f(x)=f^n(\langle x,e_1\rangle,\dots,\langle x,e_n\rangle)$, where $f^n\colon \R^n\to\R$ is the solution of the equation (\ref{plotn3}).\\

3. Suppose  $g\in D, C\in D, B\in D_H$. Let us prove that the operator  $L$ is dissipative. Let $f\in D$ and $\lambda >0$ be fixed.

As in the proof of item 2 of this theorem, the value of the function $Lf$ at the point  $x\in H$ is equal to the value of the function $L^nf^n$
at the point $(\langle x,e_1\rangle,\dots,\langle x,e_n\rangle)\in \R^n$.
Again, the operator $A$ is positive, and the function $g$ satisfies the inequality $g(x)\geq g_0\equiv\mathrm{const}>0$,
so we can apply item 2 of lemma \ref{lm_urchp} to the finite-dimensional operator. This gives us that the
finite-dimansional operator is dissipative, which implies that the operator $L$ is dissipative too. This idea can be expressed in a more formal way:

$$\left\|Lf - \lambda f\right\|=\sup_{x\in H}|g(x)\tr(Af''(x))+\langle f'(x),AB(x)\rangle +C(x)f(x) - \lambda f(x)|=$$

$$\sup_{(x_1,\dots,x_n)\in\R^n} \Bigg|g^n(x_1,\dots,x_n)\sum_{s=1}^n\sum_{k=1}^n \langle Ae_s,e_k\rangle \partial_k\partial_s f^n(x_1,\dots,x_n) +$$

$$\sum_{s=1}^n \left(\sum_{k=1}^n B^n_k(x_1,\dots,x_n)\left<Ae_s,  e_k\right>\right)   \partial_sf^n(x_1,\dots,x_n)+$$

$$C^n(x_1,\dots,x_n)f^n(x_1,\dots,x_n) - \lambda f_n(x_1,\dots,x_n)\Bigg|     \stackrel{(\ref{diss_oc})}{\geq}$$

$$\lambda\sup_{(x_1,\dots,x_n)\in\R^n} |f^n(x_1,\dots,x_n)|=\lambda\sup_{x\in H}|f(x)|=\lambda\|f\|,$$

where $x_j=\langle x,e_j\rangle$. The inequality

\be\label{diss_j}\left\|L f - \lambda f\right\|\geq \lambda\|f\|\ee

means that the operator  $L$ is dissipative.

Finally, according to proposition \ref{pr_dissop}, the closability of $L$ follows from the fact $L$ is dissipative and densely defined.\\

4. i) Suppose $g\in X$ and $C\in X$, i.e. there exists sequences $(g_j)\subset D$ and $(C_j)\subset D$ such
that $\|g- g_j\|\to 0$ and $\|C- C_j\|\to 0$. Suppose that $B(x)=0$ for each $x\in H$. The images of the operators $(\lambda I-L)$ and $(L-\lambda I)$ are equal.
As $D$ is dense in $X$, it is enough to show that the image of the operator $(L-\lambda I)$ is dense in $D$, then it will be dense in $X$. 
Let the number  $\lambda>0$ and function $\psi\in D$ be fixed. We will approximate $f$ by the values of the operator $(L-\lambda I)$.

Since $g_j\to g$ and, for every $x\in H$, the estimate $g(x)\geq g_0$ holds, it follows that there is a number $j_0'$ such that for every  $x\in H$ and all 
$j>j_0'$ we have
\be\label{oc_g_j}g_j(x)\geq \frac{g_0}{2}.\ee

We denote the operator $L$ corresponding to the functions $g_j$ and $C_j$ by the symbol $L_j$. Note that $C_j(x)\leq 0$ and $g_j(x)\geq \frac{g_0}{2}$, 
so we can apply item 2 of this theorem to the operator $L_j$. Item 2 says that image of the operator $(L_j-\lambda I)$ is equal to $D$, so
for every $j>j_0'$ there exists a function $f_j\in D$ such that

\be \label{uravn_j}L_jf_j - \lambda f_j=\psi.\ee

The goal is to prove that $L f_j-\lambda f_j\to \psi$ as $j\to\infty$. This will imply that the image of the operator $(\lambda I-L)$ is dense in $D$.

ii) Let us prepare several estimates. First, since $C_j\to C$, there is a number $j_0$ such that for every $j>j_0$

\be\label{norm_c_j}\|C_j\|\leq 2\|C\|. \ee

Second, it follows from (\ref{uravn_j}) and (\ref{diss_j}) that for every $j>j_0'$ we have

$$\lambda\|f_j\|\stackrel{(\ref{diss_j})}{\leq}\|L_jf_j - \lambda f_j\|\stackrel{(\ref{uravn_j})}{=}\|\psi\|,$$

i.e. for every $j>j_0'$

\be\label{norm_f_j} \|f_j\|\leq\frac{\|\psi\|}{\lambda}\ee.

Finally, expressing the term $\tr(Af_j'')$ by the use of the equation 
$$\psi=L_jf_j - \lambda f_j=g_j\tr(Af_j'')+(C_j-\lambda)f_j,$$ 
we find that for every $j>\max(j_0,j_0')$

\be\label{f_j_dva_shtr}\|\tr(Af_j''\|=\left\|\frac{\psi + (\lambda - C_j)f_j}{g_j}  \right\|\stackrel{(\ref{norm_f_j}), 
(\ref{norm_c_j}), (\ref{oc_g_j})}{\leq}\frac{\|\psi\|+(\lambda+ 2\|C\|)\|\psi\|/\lambda}{g_0/2}.\ee

iii) Now let us prove that if $j\to\infty,$ then $L f_j-\lambda f_j\to \psi$. Indeed, for every $j>\max(j_0,j_0')$ one obtains

$$\|L f_j-\lambda f_j- \psi\|\stackrel{(\ref{uravn_j})}{=}\|L f_j-\lambda f_j - (L_jf_j - \lambda f_j)\|=\|(g-g_j)\tr(Af_j'') + (C-C_j)f_j\|$$

$$\stackrel{(\ref{f_j_dva_shtr}), (\ref{norm_f_j})}{\leq} \|(g-g_j)\|\frac{\|\psi\|+(\lambda+ 2\|C\|)\|\psi\|/\lambda}{g_0/2} + \|(C-C_j)\|\frac{\|\psi\|}{\lambda}\to 0,$$

because $\|(g-g_j)\|\to 0$ and $\|(C-C_j)\|\to 0$. Item 4 is proven.

5. Let the coefficients of the operator $L$ be uniform limits $g, B, C$ of the continuously differentiable cylindrical functions $g_j, B_j, C_j$. 
As $g_j\to g$ and, for all $x\in H$, we have $g(x)\geq g_0$, it follows that there exists a number $j_0$ such that for all $x\in H$ and all $j>j_0$ 
we have $g_j(x)\geq \frac{g_0}{2}$. Also recall that $C_j(x)\leq 0$. This all allows us to use item 3 of this theorem.

According to (\ref{diss_j}), for every function $\varphi\in D$ and every $\lambda>0$ we have
$$\|g_j\tr(A\varphi'') + \left<\varphi',AB_j\right> + C_j\varphi -\lambda \varphi\|\geq \lambda\|\varphi\|.$$

Taking the limit as $j\to\infty$, we obtain the estimate $\|L\varphi - \lambda\varphi\| \geq \lambda\|\varphi\|,$ which means that $L$ is dissipative. 
According to proposition \ref{pr_dissop}, the dissipative operator $(L,D)$
with the domain $D$ dense in $X$ is closable. Let us denote its closure with $(\overline{L},D_1)$. Note that by proposition \ref{pr_dissop} 
the closure also is a dissipative operator. \QED

The constructions above were built to prove the following result:

\thm\label{svaz_polugr} (On the connection between the family $(S_t)_{t\geq 0}$ and the semigroup with the generator $\overline{L}$)

\textbf{Suppose} that $g\in X$, $B\in X_H$, $C\in X$, and for every $x\in H$ we have $g(x)\geq g_0\equiv\mathrm{const}>0$ and $C(x)\leq 0$. As $C\in X$, 
there exists a sequence $(C_j)\subset D$, converging to $C$ uniformly; let us additionally claim that this sequence can be selected in such a way 
that $C_j(x)\leq 0$ for all $j\in\N$ and all $x\in H$. \textbf{Then} the following holds:

1. If the closure $(\overline{L},D_1)$ of the operator $(L,D)$ is a generator of a strongly continuous semigroup $\left(e^{t\overline{L}}\right)_{t\geq 0}$ 
of linear continuous operators on the space $X$, then
\be e^{t\overline{L}}\varphi=\lim_{n\to\infty}\left(S_{\frac{t}{n}}\right)^n\varphi,\ee
where limit exists for every $\varphi\in X$ and is uniform with respect to $t\in [0,t_0]$ for every $t_0>0$.

2. If $B=0$, then the operator $(\overline{L},D_1)$ is a generator of a strongly continuous semigroup $\left(e^{\overline{L}t}\right)_{t\geq 0}$ of 
linear continuous operators on the space $X$. Moreover for every $t\geq 0$ we have 
$\left\|e^{\overline{L}t}\right\|\leq 1$, i.e. the semigroup $\left(e^{\overline{L}t}\right)_{t\geq 0}$ is contractive.

3. Suppose $B=0$, and for all $j\in\N$ the functions $g_j\in X$, $B_j\in X_H$ and $C_j\in X$ are given. Suppose $B_j=0$ for all $j\in\N$. 
Suppose there exists a number $\varepsilon_0>0$ such that for all $j\in\N$ and all $x\in H$ we have $g_j(x)\geq \varepsilon_0$ and $C_j(x)\leq 0$. 
Let us denote by the symbol $L_j$ the operator $L$, which corresponds to the functions $g_j$, $B_j$ and $C_j$, and the operator $L$ corresponding to
the functions $g$, $B$ and $C$ will be denoted by $L_0$. Suppose also that $g_j(x)\to g(x)$ and $C_j(x)\to C(x)$, uniformly with respect to  $x\in H$.

Then the (existing by item 2) strongly continuous semigroups $\left(e^{\overline{L_j}t}\right)_{t\geq 0}$ converge strongly
(and uniformly with respect to $t\in[0,t_0]$ for every fixed $t_0>0$)
to the (existing by item 2) strongly continuous semigroup $\left(e^{\overline{L_0}t}\right)_{t\geq 0}$ with the generator $\overline{L_0}$.
In other words for every $t_0>0$ and every $\varphi\in X$ there exists a limit 
\be\lim_{j\to\infty}\left(e^{\overline{L_j}t}\varphi\right)(x)=\left(e^{\overline{L_0}t}\varphi\right)(x),\ee
uniformly with respect to $x\in H$ and $t\in[0,t_0]$.

\textbf{Proof.}

1. Recall theorem \ref{FormulaChernova} and set $F(t)=S_t$, $\omega=\frac{2\|A\|\|B\|^2}{g_0}+\|C\|,$ $\mathcal{X}=X$, $\mathcal{D}=D,$ $F'(0)=L$, $G=\overline{L}$. One can see that according to items 1, 4 and 5 of 
theorem \ref{preChern2} and item 5 of theorem \ref{gDeltaA_thm} all the conditions of theorem \ref{FormulaChernova} are fulfilled.

2. Note that $C(x)\leq 0,$ so $\sup_{x\in H}e^{C(x)}\leq 1$ and for $B=0$ one obtains the estimate $\|S_t\|\leq 1$. Conditions of theorem \ref{FormulaChernova2} are fulfilled if one sets $\mathcal{X}=X$, $\mathcal{D}=D$, $\mathcal{L}=L$, $V_t=S_t$, $M=1$, $\omega=0.$ Indeed, according to item 1 of theorem \ref{preChern2}, for all $t\geq 0$ the estimate $\|S_t\|\leq e^{\omega t}=1$ holds true, therefore 
$\|\left(S_t\right)^k\|\leq 1\cdot\dots\cdot 1 =1$. Other conditions of theorem \ref{FormulaChernova2} follow from item 4 of theorem \ref{preChern2} and items 4 and 5 of theorem \ref{gDeltaA_thm}. 

3. Recall theorem \ref{TC}, and set $\mathcal{X}=X$, $\mathcal{D}=D$, $\mathcal{L}=L_0$, $\mathcal{L}_n=L_j$. One can see that 
item 2 of this theorem and items 4 and 5 of theorem \ref{gDeltaA_thm} imply all the conditions of theorem \ref{TC}, except for the following 
one: if $\varphi\in D$, then $\lim\limits_{j\to\infty}L_j\varphi=L_0\varphi$. A simple check shows that this condition is also fulfilled.
\QED

\subsection{Feynman formula solves the Cauchy problem for the parabolic equation}

We want to find a function $u\colon [0,+\infty)\times H\to\R$ satisfying the following conditions (we call them 
Cauchy problem for the parabolic differential equation):

\be\label{CP}
\left\{ \begin{array}{ll}
 u'_t(t,x)=Lu(x,t);& \quad  t\geq 0, x\in H,\\
 u(0,x)=u_0(x); & \quad  x\in H.\\
  \end{array} \right.
\ee

To this Cauchy problem, we relate the so-called abstract Cauchy problem (see Definition \ref{abszadC}), which we define as the following system of conditions upon the function $U\colon [0,+\infty)\to X$:
\be\label{ACP}
\left\{ \begin{array}{ll}
 \frac{d}{d t}U(t)= \overline{L}U(t); & t\geq 0,\\
 U(0)=u_0,\\
  \end{array} \right.
\ee

\rem\label{zam_sved_cp} Problem (\ref{CP}) can be considered as problem (\ref{ACP}) in the following sense. 
Function $u\colon (t,x)\longmapsto u(t,x)$ of two variables $(t,x)$ can be considered as a function $u\colon t\longmapsto [x\longmapsto u(t,x)]$ of one variable $t$,
with values in the space of functions of variable $x$.
Then
$$u(t,x)=(U(t))(x),\quad t\geq 0, x\in H.$$
Using this  correspondence, we start from Definition \ref{abszadC} and define the solution of problem (\ref{CP}).

\df We call a function $u\colon [0,+\infty)\times H\to \R$ \textit{a strong solution} of problem (\ref{CP}) if it satisfies the following conditions:

\be\label{CPreg}
\left\{ \begin{array}{ll}
 u(t,\cdot )\in D_1; & t\geq 0,\\
\textrm{function } t\longmapsto u(t,\cdot) \textrm{ is continuous};& t\geq 0,\\
 \textrm{Uniformly for } x\in H\ \exists \lim\limits_{\varepsilon\to 0}\frac{u(t+\varepsilon,x)-u(t,x)}{\varepsilon}=u'_t(t,x); & t\geq 0,\\

u'_t(t,\cdot)\in X;& t\geq 0,\\
\textrm{Function } t\longmapsto u'_t(t,\cdot)\textrm{ is continuous};& t\geq 0,\\
 u'_t(t,x)=Lu(x,t); & t\geq 0, x\in H,\\
 u(0,x)=u_0(x); & x\in H.\\
  \end{array} \right.
\ee

\df We call a function $u\colon [0,+\infty)\times H\to \R$ \textit{a mild solution} of problem (\ref{CP}) if it satisfies the following conditions:

\be\label{CPmild}
\left\{ \begin{array}{ll}
u(t,\cdot)\in X; & t\geq 0,\\
\textrm{Function } t\longmapsto u(t,\cdot) \textrm{ is continuous};& t\geq 0,\\
 \int_0^t u(s,\cdot)ds \in D_1; & t\geq 0,\\
 u(t,x)=L\int_0^tu(s,x)ds + u_0(x); & t\geq 0, x\in H,\\
 u_0\in X.\\
  \end{array} \right.
\ee

\df\label{klass} Let us use the symbol $C([0,+\infty),X)$ for the class of all functions $u\colon [0,+\infty)\times H\to\R$ such that for every $t\geq 0$ 
the function $x\longmapsto u(t,x)$ belongs to the class $X$, and the mapping $t\longmapsto u(t,\cdot)\in X$ is continuous for every $t\geq 0$.\\

Finally, let us state and prove the main result of the article. We use definitions and notation from Section \ref{denotions}.

\thm\label{main_th_art}  (On the solution of the Cauchy problem for a parabolic differential equation in Hilbert space)

Suppose $g\in X, C\in X, B\in X_H$. Suppose there is a number $g_0>0$ such that for all $x\in H$ we have $g(x)\geq g_0$ and $C(x)\leq 0$. As $C\in X$, 
there exists a sequence $(C_j)\subset D$, converging to $C$ uniformly; let us additionally require
that this sequence can be selected in such a way way that $C_j(x)\leq 0$ for all $j\in\N$ and all $x\in H$.

Then the following holds:

1. If there exists a strongly continuous semigroup with the generator $\overline{L}$, then for every $u_0\in D_1$ there exists a solution $u$ 
of problem (\ref{CPreg}), unique in the class $C([0,+\infty),X)$. The solution depends continuously on $u_0$, and is given by
the formula $u(t,x)=\lim\limits_{n\to\infty}\left(\left(S_{\frac{t}{n}}\right)^nu_0\right)(x),$ where the limit is uniform with respect 
to $t\in [0,t_0]$ for every $t_0>0$.

2. If there exists a strongly continuous semigroup with the generator $\overline{L}$, then for every $u_0\in X$ there exists a solution $u$ of 
problem (\ref{CPmild}), unique in the class $C([0,+\infty),X)$. It depends continuously on $u_0$, and is  given by
the formula $u(t,x)=\lim\limits_{n\to\infty}\left(\left(S_{\frac{t}{n}}\right)^nu_0\right)(x),$ where the limit is uniform with respect to $t\in [0,t_0]$ 
for every $t_0>0$.

3. If $B=0$, then there exists a strongly continuous semigroup with the generator $\overline{L}$. 
The formula $u(t,x)=\lim\limits_{n\to\infty}\left(\left(S_{\frac{t}{n}}\right)^nu_0\right)(x)$ becomes simpler than in the case $B\neq 0$. Namely, for $B=0$ we have

\be\label{FF1}u(t,x)=\lim_{n\to\infty} \underbrace{\int_H\int_H\dots\int_H\int_H}_ne^{\frac{t}{n}\left( C(x)+\sum_{k=1}^{n-1}C(y_k)\right)} 
u_0(y_1)\mu_{\frac{2t}{n}g(y_2)A}^{y_2}(dy_1)\mu_{\frac{2t}{n}g(y_3)A}^{y_3}(dy_2)\dots\ee
$$ \dots\mu_{\frac{2t}{n}g(y_n)A}^{y_n}(dy_{n-1}) \mu_{\frac{2t}{n}g(x)A}^{x}(dy_n).$$

In this case the solution $u$ for all $t>0$ satisfies the estimate $\sup_{x\in H}|u(t,x)|\leq \sup_{x\in H}|u_0(x)|$.

4. Let $B=0$, and let the functions  $g_j\in X$, $B_j\in X_H$ and $C_j\in X$ be given for all $j\in\N$. Let $B_j=0$ for all $j\in\N$. 
Suppose there exists $\varepsilon_0>0$ such that $g_j(x)\geq \varepsilon_0$ and $C_j(x)\leq 0$ for all $j\in\N$ and all $x\in H$.
Let us use the symbol $L_j$ for the operator $L$ that corresponds to the functions $g_j$, $B_j$ and $C_j$, and the symbol $L_0$ for the operator $L$
that corresponds to the functions $g$, $B$ and $C$. Suppose also that $g_j(x)\to g(x)$ and $C_j(x)\to C(x)$, uniformly with respect to $x\in H$. 
We denote as $u_j$ the solution of problems (\ref{CPreg}) and (\ref{CPmild}) for the operator $L_j$. For solution of problems (\ref{CPreg}) and (\ref{CPmild})
with the operator $L$, we use the symbol $u$.

Then $u_j(t,x)$ converges to $u(t,x)$ as $j\to\infty$, uniformly with respect $x\in H$ and uniformly with respect to $t\in[0,t_0]$ for every fixed $t_0>0$.

\rem Note that if $B=0$, then solution depends continuously on the data of the Cauchy problem: the coefficients of the equation (item 4) and 
the initial condition (items 1 and 2).

\rem Analogous theorems for $\mathbb{C}$- or $\R^n$-valued functions $u$
can be formulated mutatis mutandis. The result will hold true due to the theorem above and the linearity of $L$ and $S_t$. 
The only additional condition will be that the coefficients of the equation must be real-valued. The same remark is applicable 
to all the key theorems of this article.\\

\textbf{Proof of the theorem.}

1. Suppose that there exists a strongly continuous semigroup with the generator $\overline{L}$. Then by item 1 of proposition \ref{ACPsol} 
we obtain the existence of a strong solution (definition \ref{abszadC}) to Cauchy problem (\ref{ACP}),  and the solution is unique in the class $C([0,+\infty), X)$. 
By item 1 of theorem \ref{svaz_polugr} the semigroup is given in the form described. Using the relation between problems (\ref{CP}) and (\ref{ACP}) 
explained in remark \ref{zam_sved_cp}, we obtain the solution for problem (\ref{CPreg}). The solution is unique in the class $C([0,+\infty), X)$, as follows 
from remark \ref{zam_sved_cp}.

2. The proof is similar to that in item 1. The only difference is that in proposition \ref{ACPsol} we use item 2 instead of item 1.

3. The existence of the sought semigroup follows from item 2 of theorem \ref{svaz_polugr}. The estimate for the supremum of the absolute value of the solution 
follows from the fact that the semigroup is contractive.

Let us explain how the equality $u(t,x)=\lim_{n\to\infty} \left(\left(S_{\frac{t}{n}}\right)^nu_0\right)(x)$ implies formula (\ref{FF1}).
For a continuous bounded function $\psi\colon H\to\R$ and a point $x\in H$, the following change of variables rule in the integral is correct:
$$\int_H\psi(y)\mu_A(dy)=\int_H\psi(y-x)\mu_A^{x}(dy).$$
Applying this rule, and changing $A$ to ${2tg(x)A}$, we come to the equality
$$(S_t\varphi)(x)=e^{tC(x)}\int_H\varphi(x+y)\mu_{2tg(x)A}(dy)$$ $$=e^{tC(x)}\int_H\varphi(x+(y-x))\mu_{2tg(x)A}^x(dy)=e^{tC(x)}\int_H\varphi(y)\mu_{2tg(x)A}^{x}(dy).$$
For $n=2$ in formula (\ref{FF1}) we get the expression
$$\left(\left(S_{\frac{t}{2}}\right)^2\varphi\right)(x)=\left(S_{\frac{t}{2}}\left(S_{\frac{t}{2}}\varphi\right)\right)(x)
=\int_H\left( \int_He^{\frac{t}{2}\left( C(x)+C(y_1)\right)}\varphi(y_1)\mu_{\frac{2t}{2}g(y_2)A}^{y_2}(dy_1)\right)\mu_{\frac{2t}{2}g(x)A}^{x}(dy_2).$$
In the same way expressions for  $n>2$ are derived. Thus, the formula (\ref{FF1}) is proven.

4. The proof follows immediately from item 3 of theorem \ref{svaz_polugr}.
\QED

\section{Acknowledgments}

The author is grateful to O.G.~Smolyanov for setting the problem (in particular, for defining the space $X$) 
and attention to the work; to V.I.~Bogachev, O.G.~Smolyanov and E.T.~Shavgulidze for help with the calculation of integral (\ref{int_eksp2}); 
to T.A.~Shaposhnikova for acquainting the author with book \cite{Kryl}, where the key lemma \ref{lm_urchp} was taken from; to A.V.~Halyavin 
for reading and commenting on the manuscript.

The author also highly appreciates the contributions of L.S.~Remizova, Yu.A.~Komlev, D.V.~Turaev, A.S.~Remizova and  A.G.~Linucheva,
whose kind support was a necessary condition of the present paper's existence.

This work has been supported by the Russian Scientific Foundation Grant 14-41-00044 at the Lobachevsky University of Nizhny Novgorod.


\begin{thebibliography}{199}

\bibitem{Bog2} Bogachev V.I. Gaussian Measures. --- Amer. Math. Soc., Providence, 1998.


\bibitem{BS}
V.I. Bogachev, O.G. Smolyanov. Real and functional analysis: university course. --- M. Izevsk: RCD, 2009.

\bibitem{Butko1} Ya.A. Butko. The Feynman-Kac-Ito formula for an infinite-dimensional Schr\"{o}dinger equation with scalar and vector potentials.

\bibitem{Butko-mult} Ya.A. Butko. Feynman formula for semigroups with multiplicatevely perturbed generators. --- Science and Education (ISSN 1994-0408), $\#$10, October 2011. Article number 77-305691/239563.

\bibitem{DF} Yu.L. Daletsky, S.V. Fomin. Measures and differential equations in infinite-dimensional space. --- Kluwer, 1991.

\bibitem{EgZhLob} A.D. Egorov, E.P. Zidkov, Yu.Yu. Lobanov. Introduction to the theory and applications of the functional integration (in Russian). --- M. Fizmatlit, 2006.

\bibitem{Kryl} N.V. Krylov. Lectures on Elliptic and Parabolic Equations in Holder Spaces. --- AMS, Graduate Texts in Mathematics vol. 12, 1996.

\bibitem{Lobanov} Yu.Yu. Lobanov. Methods of the approximate functional integrating for the numerical research in the quantum physics. Dr. Sci. dissertation thesis (in Russian). --- M.,2009.


\bibitem{Sm1} O.G. Smolyanov. Analysis on the topological linear spaces and its applications (in Russian). --- M.: MSU, 1979.

\bibitem{SSh} O.G. Smolyanov, E.T. Shavgulidze. Continual integrals (in Russian). --- M.: MSU, 1990.

\bibitem{SShK} Smolyanov O.G., Shamarov N.N., Kpekpassi M. Feynman-Kac and Feynman Formulas for Infinite-Dimensional Equations with Vladimirov Operator. --- Doklady Mathematics, 83:3 (2011), 389–393.

\bibitem{SSham} O.G. Smolyanov and N.N. Shamarov. Hamiltonian Feynman formulas for equations containing the Vladimirov operator with variable coefficients. --- Doklady Mathematics, Volume 84, Number 2, 689-694, DOI: 10.1134/S1064562411060330

\bibitem{Schwartz}  L. Schwartz. Analyse mathematique, I. Hermann, 1967.

\bibitem{Botelho} Luiz C.L. Botelho. Non-linear diffusion in $\R^D$ and Hilbert Spaces, a Cylindrical/Functional Integral Study. --- arXiv:1003.0048v1 [physics.gen-ph] 27 Feb 2010.


\bibitem{BGS2010} Ya.A. Butko, M. Grothaus, O.G. Smolyanov. Lagrangian Feynman formulas for second-order parabolic equations in bounded and unbounded domains. --- Infinite Dimansional Analyasis, Quantum Probability and Related Topics, vol. 13, No. 3 (2010), 377-392.

\bibitem{BShS} Yana A. Butko, Ren\'e L. Schilling, Oleg G. Smolyanov. Lagrangian and Hamiltonian Feynman formulae for some Feller semigroups and their perturbations. --- arXiv:1203.1199v1 [math.PR] 6 Mar 2012.

\bibitem{Cartan} H. Cartan. Differential Calculus. --- Kershaw Publishing Company, 1971.

\bibitem{Pazy} A. Pazy. Semigroups of Linear Operators and Applications to Partial Differential Equations. --- Springer-Verlag, 1983.

\bibitem{Plya1} A.S. Plyashechnik. Feynman formula for Schredinger-Type equations with time- and space-dependent coefficients,
Russian Journal of Mathematical Physics, 2012, vol. 19, No.3, pp. 340-359.

\bibitem{Plya2} A.S. Plyashechnik. Feynman formulas for second-order parabolic equations with variable coefficients,
Russian Journal of Mathematical Physics, 2013, vol. 20, No.3, pp. 377-379.

\bibitem{DaPrato} G. Da Prato. Introduction to infinite-dimensional analysis. --- Springer, 2006.

\bibitem{DaPrato2} G. Da Prato, J. Zabczyk. Second Order Partial Differential Equations in Hilbert Spaces. --- London Mathematical Society Lecture Notes Series 293, 2004.



\bibitem{Fern} X. Fernique. Int\'egrabilit\'e des vecteurs gaussiens. --- C. R. Acad. Sci. Paris S\'er. A-B 270: A1698–A1699. See also MR266263 (42 \#1170) 60.08   (1970).

\bibitem{F1} R.P. Feynman. Space-time approach to nonrelativistic quantum mechanics. --- Rev. Mod. Phys., 20 (1948), 367-387.

\bibitem{F2} R.P. Feynman. An operation calculus having applications in quantum electrodynamics. --- Phys. Rev. 84 (1951), 108-128.


\bibitem{EN} K.-J. Engel, R. Nagel. One-Parameter Semigroups for Linear Evolution Equations. --- Springer, 2000.

\bibitem{Djedone} J. Dieudonn\'e. Foundations of modern analysis. --- Academic Press, New York and London, 1969.


\bibitem{Remizov} I.D. Remizov. Solution of a Cauchy problem for a diffusion equation in a Hilbert space by a Feynman formula// Russian Journal of Mathematical Physics, 2012, v.19, No.3, 360-372.


\bibitem{Simon} B. Simon. Functional Integration and Quantum Physics. --- Academic Press, 1979.

\bibitem{Sm2} O.G. Smolyanov, A.G. Tokarev, A. Truman. Hamiltonian Feynman path integrals via the Chernoff formula. --- J. Math. Phys. 43, 10 (2002) 5161-5171.

\bibitem{Sm3} O.G. Smolyanov. Feynman formulae for evolutionary equations. --- Trends in Stochastic Analyasis, London Mathematical Society Lecture Notes Series 353, 2009.

\bibitem{Go} H.-S. Kuo. Gaussian measures in Banach space. --- Lecture notes in mathematics, 463. Springer-Verlag, 1975.

\bibitem{Chernoff} Paul R. Chernoff, Note on product formulas for operator semigroups, J. Functional Analysis 2 (1968), 238-242. 

\bibitem{SWW}O.G. Smolyanov, H. v. Weizs\"{a}cker, O. Wittich. Chernoff's Theorem and Discrete Time Approximations of Brownian Motion on Manifolds//Potential Analysis, February 2007, Volume 26, Issue 1, pp 1-29.

\bibitem{Butko4} Ya.A. Butko. Feynman formulae for evolution semigroups  (in Russian). Electronic scientific and technical periodical "Science and education", DOI: 10.7463/0314.0701581 , N 3 (2014), 95-132. 

\bibitem{OSS} Yu.N. Orlov, V.Zh. Sakbaev, O.G. Smolyanov. Feynman formulas as a method of averaging random Hamiltonians.//Proceedings of the Steklov Institute of Mathematics, August 2014, Volume 285, Issue 1, pp 222-232.

\bibitem{Rmz1} I.D. Remizov. The latest version of the preprint http://arxiv.org/abs/1409.8345

\bibitem{Bot1} L.C.L. Botelho. Non linear Diffusion and Wave Dumped Propagation: Weak Solutions and Statistical Turbulence Behavior.// Journal of Advanced Mathematics and Applications, vol. 3, 1-11, 2014.

\bibitem{Bot2} L.C.L. Botelho. Semi-linear Diffusion in $R^D$ and Hilbert Spaces, a Feynman-Wiener path integral study. --- Random Oper. Stoch. Equ., doi 10.1515/ROSE.2011.020 (2011)

\bibitem{Bot3} L.C.L. Botelho. A method of integration for wave equation and some applications to wave physics. --- Random Oper. Stoch. Equ., 18 (2010), pp. 301-325, doi 10.1515/ROSE.2010.017.

\bibitem{Bot4} L.C.L. Botelho. A note on Feynman-Kac path integral representations for scalar wave motions. --- Random Oper. Stoch. Equ., doi 10.1515/rose-2013-0012 (2013)

\end{thebibliography}
\end{document}